\documentclass[11pt,openbib]{article}

\usepackage{amsmath,amssymb,theorem}
\usepackage[mathscr]{eucal}

\newtheorem{theorem}{\sc Theorem}[section]
\newtheorem{lemma}[theorem]{\sc Lemma}
\newtheorem{proposition}[theorem]{\sc Proposition}
\newtheorem{corollary}[theorem]{\sc Corollary}

{\theorembodyfont{\rmfamily}

\newtheorem{example}[theorem]{\sc Example}
\newtheorem{remark}[theorem]{\sc Remark}}

\makeatletter \@addtoreset{equation}{section} \makeatother

\newcommand{\book}[1]{}

\title{\huge\bf
A QUANTUM CANONICAL EMBEDDING}
\author{\Large\bf Leonid L. Vaksman}
\date{}

\begin{document}


\maketitle

\begin{center}
\Large with the participation of
\\ O. A. Bershtein, Ye. K. Kolisnyk, and S. D. Sinel'shchikov
\end{center}



\setcounter{section}{-1}

\section{Foreword}

These notes constitute a compilation of a part of a book by Leonid Vaksman `Quantum
Bounded Symmetric Domains'. The book was in the process of writing down during the last
decade and is concentrated on expounding the results of L. Vaksman and his team on
quantum bounded symmetric domains, along with a large variety of mathematical background
and preliminary information. The Russian text of this book can be found at
arXiv:0803.3769 [math.QA]. It is frustrating to observe that the above monograph has not
been finished for publication purposes before L. Vaksman passed away so suddenly.
Nevertheless, a great deal of the material is available in preprints and journal papers.
Perhaps, the only essential exception here is the chapter on the canonical embedding (for
a special case we refer to \cite{SSV4}), and so the principal motive of the people listed
as participants was to fill this gap.

The canonical embedding we are about to produce is a q-analog for a map of
complex varieties which derives the Borel embedding of a bounded symmetric
domain from the Harish-Chandra embedding. Additionally, we produce a
quantum analog of the Harish-Chandra embedding of the associated real
varieties, see \eqref{classic_embedds} below.

The exposition is preceded with some background material as well as the simple example of
the quantum disc. This section has been added recently for the convenience of the reader.

Special thanks are to be expressed to Prof. H. P. Jakobsen for many
fruitful discussions on the entire text.

\section{Preliminaries}

In this section we list some important results and notations from the theory of quantum
bounded symmetric domains. All these facts can be found in various papers, but we list
them here for the sake of convenience.

Let $(a_{ij})_{i,j=1,\ldots,l}$ be a Cartan matrix and $\mathfrak g$ the corresponding
simple complex Lie algebra. Hence the Lie algebra can be defined by the generators
$e_i,f_i,h_i,i=1,...,l$, and the well-known relations, see \cite{Jant}. Let
$\mathfrak{h}$ be the linear span of $h_i,i=1,...,l$. The simple roots $\{\alpha_i \in
\mathfrak{h}^*|i=1,...,l\}$ are given by $\alpha_i(h_j)=a_{ji}$. Also, let
$\{\overline{\omega}_i| i=1,...,l\}$ be the fundamental weights, hence
$P=\bigoplus_{i=1}^l
\mathbb{Z}\overline{\omega}_i=\{\lambda=(\lambda_1,...,\lambda_l)|\,\lambda_j \in
\mathbb{Z}\}$ is the weight lattice and $P_+=\bigoplus_{i=1}^l
\mathbb{Z_+}\overline{\omega}_i=\{\lambda=(\lambda_1,...,\lambda_l)|\,\lambda_j \in
\mathbb{Z}_+\}$ is the set of integral dominant weights.

Fix $l_0 \in \{1,...,l\}$, $\mathbb S=\{1,...,l\} \backslash \{l_0\}$, and the Lie
subalgebra $\mathfrak{k} \subset \mathfrak{g}$ generated by
$$
e_i, f_i,\;\; i\neq l_0; \qquad h_i, \;\;i=1,...,l.
$$
Define $H_{\mathbb S} \in \mathfrak{h}$ by
$$
\alpha_i(H_{\mathbb S})=0, \; i \neq l_0; \qquad
\alpha_{l_0}(H_{\mathbb S})=2.
$$
We restrict ourselves to Lie algebras $\mathfrak{g}$ that can be
equipped with a $\mathbb Z$-grading as follows:
\begin{equation}\label{par_type}
\mathfrak{g}=\mathfrak{g}_{-1}\oplus\mathfrak{g}_0\oplus\mathfrak{g}_{+1},
\quad \mathfrak{g}_j=\{\xi\in\mathfrak{g}|\:[H_{\mathbb
S},\xi]=2j\xi\}.
\end{equation}
Let $\delta$ be the maximal root, and $\delta=\sum_{i=1}^l c_i\alpha_i$. In this setting
(\ref{par_type}) holds if and only if $c_{l_0}=1$. In this case
$\mathfrak{k}=\mathfrak{g}_0$, and we use more conventional notation
$\frak{p}^{\pm}\stackrel{\rm def}{=}\frak{g}_{\pm 1}$. The pair
$(\mathfrak{g},\mathfrak{k})$ is called a Hermitian symmetric pair.

Recall two important embeddings that arise in the theory of bounded symmetric domains
(see \cite{Wolf}, \cite{Kor}). The Harish-Chandra embedding gives us the canonical
realization of a bounded symmetric domain $\mathbb D$ as a unit ball in $\mathfrak{p}^-$.
On the other hand, $\mathbb D =G_{\mathbb R }/K_{\mathbb R}$ is a Hermitian symmetric
space of noncompact type, where $G_{\mathbb R}$ is the group of isometries of $\mathbb D$
and $K_{\mathbb R}$ is the stabilizer of a point. Let $P \subset G$ be the connected
complex Lie subgroup with $\operatorname{Lie}P=\mathfrak{k}\oplus \mathfrak{p}^+$. The
Borel embedding theorem claims that $\mathbb D$ can be embedded into the projective
variety $X=G/P$ as an open $G_{\mathbb R}$-orbit. Moreover, we have the following
embeddings:
\begin{equation}\label{classic_embedds}
G_{\mathbb R}/K_{\mathbb R} = \mathbb D \hookrightarrow \mathfrak{p}^- \hookrightarrow
G/P.
\end{equation}
In this paper we develop a quantum analog for the second embedding.

\subsection{\boldmath An algebra $\mathrm{Pol}(\mathfrak{p}^-)_q$}

Now we recall some basic notions from the theory of quantum bounded symmetric domain. We
describe the construction of a noncommutative associative algebra
$\mathrm{Pol}(\mathfrak{p}^-)_q$ which is considered as a quantum analog of the algebra
of polynomials on $\mathfrak{p}^-$. This notion was introduced by three groups of
mathematicians independently. Kamita, Morita and Tanisaki established quantum analogs of
some prehomogeneous spaces of commutative parabolic type (see \cite{KMT}). Jakobsen
developed quantum Hermitian symmetric spaces in \cite{Jak-Hermit}. Finally, Vaksman and
Sinel'shchikov managed to introduce quantum analogs of bounded symmetric domains. The
equivalence of these approaches to the construction of $\mathrm{Pol}(\mathfrak{p}^-)_q$
was proved by D. Shklyarov in \cite{Shkl}. Here we recall the Vaksman-Sinel'shchikov
construction of the $U_q \mathfrak{g}$-module algebras $\mathrm{Pol}(\mathfrak{p}^-)_q$
briefly, we refer the readers to \cite{SV} for the proofs.

We fix in the sequel $q \in (0,1)$. Consider the Hopf algebra $U_q \mathfrak{g}$ over
${\Bbb C}$ given by the generators $\{E_i,F_i,K_i^{\pm 1}\}_{i=1}^l$ and the
Drinfeld-Jimbo relations
$$K_iK_j=K_jK_i,\quad K_iK_i^{-1}=K_i^{-1}K_i=1,$$
\begin{equation*}
K_iE_j=q_i^{a_{ij}}E_jK_i,\quad K_iF_j=q_i^{-a_{ij}}F_jK_i,
\end{equation*}
\begin{equation*}
E_iF_j-F_jE_i=\delta_{ij}\,\frac{K_i-K_i^{-1}}{q_i-q_i^{-1}},
\end{equation*}
\begin{equation*}
\sum\limits_{m=0}^{1-a_{ij}}(-1)^m\left[1-a_{ij}\atop m\right]_{q_i}
E_i^{1-a_{ij}-m}E_jE_i^m=0,
\end{equation*}
\begin{equation*}
  \sum\limits_{m=0}^{1-a_{ij}}(-1)^m\left[1-a_{ij} \atop m\right]_{q_i}
  F_i^{1-a_{ij}-m}F_jF_i^m=0,
\end{equation*}
where $d_i>0$, $i=1,2,\ldots,l$, are coprime integers such that the matrix
$\mathbf{b}=(d_i a_{ij})_{i,j=1,2,\ldots,l}$ is positive definite, $q_i=q^{d_i}$, $1\leq
i\leq l$,
\begin{equation*}
\left[m \atop n\right]_q=\frac{[m]_q!}{[n]_q![m-n]_q!},\quad
[n]_q!=[n]_q\ldots[2]_q[1]_q,\quad [n]_q=\frac{q^n-q^{-n}}{q-q^{-1}}.
\end{equation*}

We equip this Hopf algebra with a grading as follows:
$${\rm deg}\,K_j\:=\:{\rm deg}\,E_j\:=\:{\rm deg}\,F_j\:=\:0,\qquad j \ne
l_0$$
$${\rm deg}\,K_{l_0}\,=\,0,\quad {\rm deg}\,E_{l_0}\,=\,1,\quad {\rm
deg}\,F_{l_0}\,=\,-1.$$

Recall certain important notions. First of all, a $U_q \mathfrak{g}$-module $V$ is called
weight if it decomposes into a direct sum of subspaces as follows:
\begin{equation*}
V=\bigoplus_{\lambda \in P}V_\lambda, \qquad V_\lambda=\{v \in V|
K_jv=q_j^{\lambda_j}v,j=1,...,l\}.
\end{equation*}

For any $\lambda \in P_+$ denote by $L(\lambda)$ a $U_q \mathfrak{g}$-module with one
generator $v(\lambda)$ and the defining relations
\begin{equation*}
E_iv(\lambda)=0,\quad K_i^{\pm 1}v(\lambda)=q_i^{\pm\lambda_i}\; v(\lambda),\quad
F_i^{\lambda_i+1}v(\lambda)=0, \quad i=1,2,\ldots,l.
\end{equation*}
$L(\lambda)$ is a simple weight finite dimensional $U_q \mathfrak{g}$-module \cite{Jant}.

Recall the notion of a generalized Verma modules, see \cite{Lep}. Denote by
$U_q\mathfrak{q}^+$, $U_q\mathfrak{q}^-$ and $U_q\mathfrak{k}$ the Hopf subalgebras
generated by
$$
\{K_i^{\pm 1}, E_i\}_{i=1,2,\ldots,l}\; \cup\; \{F_j\}_{j \neq l_0}, \qquad \{K_i^{\pm
1}, F_i\}_{i=1,2,\ldots,l}\;\cup\; \{E_j\}_{j \neq l_0},
$$
$$
\{K_i^{\pm 1}\}_{i=1,2,\ldots,l}\;\cup\;\{E_j, \; F_j\}_{j \neq l_0},
$$
respectively. Evidently, $U_q\mathfrak{k}=U_q\mathfrak{q}^+\cap U_q \mathfrak{q}^-$.

Easy observations show that generalized Verma modules can be defined in terms of
generators and relations. Namely, denote by $N(\mathfrak{q}^+,0)$ a
$U_q\mathfrak{g}$-module with the generator $v(\mathfrak{q}^+,0)$ and the defining
relations
\begin{equation*}
\begin{gathered}
E_j\,v(\mathfrak{q}^+,0)=0,\quad K_j^{\pm 1}\,v(\mathfrak{q}^+,
0)=\,v(\mathfrak{q}^+,0),\qquad j=1,2,\ldots,l;
\\ F_j\,v(\mathfrak{q}^+,0)=0,\qquad j\neq l_0.
\end{gathered}
\end{equation*}
$N(\mathfrak{q}^+,0)$ is called a generalized Verma module with highest weight $0$.
Similarly, denote by $N(\mathfrak{q}^-,0)$ a $U_q\mathfrak{g}$-module with the generator
$v(\mathfrak{q}^-,0)$ and the defining relations
\begin{equation*}
\begin{gathered}
F_j\,v(\mathfrak{q}^-,0)=0,\quad K_j^{\pm 1}\,v(\mathfrak{q}^-,
0)=\,v(\mathfrak{q}^-,0),\qquad j=1,2,\ldots,l;
\\ E_j\,v(\mathfrak{q}^-,0)=0,\qquad j\neq l_0.
\end{gathered}
\end{equation*}
$N(\mathfrak{q}^-,0)$ is called a generalized Verma module with lowest weight $0$.

The $U_q \mathfrak{g}$-modules $N(\mathfrak{q}^\pm,0)$ are weight. Moreover,
\begin{equation*}
N(\mathfrak{q}^+,0)=\bigoplus_{j \in \mathbb Z_+}N(\mathfrak{q}^+,0)_{-j}, \qquad
N(\mathfrak{q}^+,0)_{-j} \:=\:\bigoplus_{\{\mu \in{\frak h}^*|\,\mu(H_{\mathbb
S})=2j\}}N(\mathfrak{q}^+,0)_\mu,
\end{equation*}
\begin{equation*}
N(\mathfrak{q}^-,0)=\bigoplus_{j \in \mathbb Z_+}N(\mathfrak{q}^-,0)_j, \qquad
N(\mathfrak{q}^-,0)_j \:=\:\bigoplus_{\{\mu \in{\frak h}^*|\,\mu(H_{\mathbb
S})=2j\}}N(\mathfrak{q}^-,0)_\mu,
\end{equation*}
with $N(\mathfrak{q}^-,0)_j$, $N(\mathfrak{q}^+,0)_{-j}$ being finite dimensional.

Let $U_q{\frak g}^{cop}$ stand for the Hopf algebra derived from $U_q{\frak g}$ by
replacing its comultiplication with the opposite one.

We intend to use the generalized Verma modules for producing coalgebras dual to covariant
algebras. To provide a precise correspondence between these two notions, we are going to
replace $U_q{\frak g}$ by $U_q{\frak g}^{cop}$ in tensor products of generalized Verma
modules.

Consider the $U_q{\frak g}$-modules $N(\mathfrak{q}^\pm,0)$. There exist evident
morphisms of $U_q{\frak g}^{cop}$-modules:
$$
\Delta^\pm:\:N(\mathfrak{q}^\pm,0)\,\to \,N(\mathfrak{q}^\pm,0)\otimes
N(\mathfrak{q}^\pm,0);\quad \varepsilon^\pm:\:N(\mathfrak{q}^\pm,0)\to \mathbb C
$$
such that
$$
\Delta^\pm: v(\mathfrak{q}^\pm, 0) \mapsto v(\mathfrak{q}^\pm, 0) \otimes
v(\mathfrak{q}^\pm, 0),\quad \varepsilon^\pm: v(\mathfrak{q}^\pm, 0) \mapsto 1.
$$

Just as in the case $q=1$, one can verify that the maps $\Delta^\pm$ are coassociative,
and that $\varepsilon^\pm$ are the counits for the coalgebras corresponding to
$\Delta^\pm$, respectively. These constructions equip $N(\mathfrak{q}^\pm,0)$ with the
structures of a $U_q{\frak g}^{cop}$-module coalgebra.

Introduce the notation
$$
\mathbb{C}[\mathfrak{p}^-]_q=\bigoplus\limits_{j\ge 0}
\mathbb{C}[\mathfrak{p}^-]_{q,j},\qquad\mathbb{C}[\mathfrak{p}^-]_{q,j}
\stackrel{\mathrm{def}}{=}N(\mathfrak{q}^+,0)_{-j}^*,
$$
$$
\mathbb{C}[\mathfrak{p}^+]_q=\bigoplus\limits_{j\ge 0}
\mathbb{C}[\mathfrak{p}^+]_{q,-j},\qquad\mathbb{C}[\mathfrak{p}^+]_{q,-j}
\stackrel{\mathrm{def}}{=}N(\mathfrak{q}^-,0)_j^*,
$$

By standard reasons, $\mathbb{C}[\mathfrak{p}^-]_q$ and $\mathbb{C}[\mathfrak{p}^+]_q$
are $U_q \mathfrak{g}$-module algebras (i.e. the multiplications are morphisms of $U_q
\mathfrak{g}$-modules, one can find the details in \cite{SV}). These covariant algebras
may be considered as q-analogs of polynomial algebras (holomorphic or antiholomorphic
identified by `-' or `+', respectively) on the quantum prehomogeneous space ${\frak
p}^-$.

\bigskip

Equip $\mathbb{C}[\mathfrak{p}^- \oplus \mathfrak{p}^+]_q \stackrel{\rm def
}{=}\mathbb{C}[\mathfrak{p}^-]_q \otimes\mathbb C[\mathfrak{p}^+]_q$ with an algebra
structure by defining a multiplication
$$
{m:\mathbb{C}[\mathfrak{p}^-]_q \oplus \mathfrak{p}^+]_q\otimes
\mathbb{C}[\mathfrak{p}^-\oplus\mathfrak{p}^+]_q\to
\mathbb{C}[\mathfrak{p}^-\oplus\mathfrak{p}^+]_q}.
$$
Namely, set
\begin{equation*}
m=(m^-\otimes m^+)(\mathrm{id}_{\mathbb{C}[\mathfrak{p}^-]_q}\otimes
\check{R}\otimes\mathrm{id}_{\mathbb{C}[\mathfrak{p}^+]_q}),
\end{equation*}

Here $m_+,\,m_-$ are the multiplications in $\mathbb C[{\frak p}^-]_q,\,\mathbb C[{\frak
p}^+]_q$ respectively, and $\check{R}:\;\mathbb C[{\frak p}^+]_q \otimes \mathbb C[{\frak
p}^-]_q \to \mathbb C[{\frak p}^-]_q \otimes \mathbb C[{\frak p}^+]_q$ is the morphism of
$U_q{\frak g}$-modules defined below by use of V. G. Drinfeld's universal R-matrix
\cite{DrinfEng}.

It is easy to show that the universal R-matrix determines a linear operator in $\mathbb
C[{\frak p}^+]_q \otimes \mathbb C[{\frak p}^-]_q$.

Now we are in a position to define the operator $\check{R}$ in a standard way:
$\check{R}=\sigma \circ R$ with $\sigma:\:a \otimes b \mapsto b \otimes a$ being the
permutation of tensors. Thus $\check{R}$ becomes a morphism of $U_q{\frak g}$-modules
since \cite{DrinfEng, Drinf2, ChP} $R$ intertwines $\Delta$ and $\Delta^{op}$.

The associativity of the multiplication in ${\rm Pol}({\frak p}^-)_q$, the existence of a
unit and the covariance of ${\rm Pol}({\frak p}^-)_q$ are proved in \cite{SV}.

\medskip

Equip $U_q \mathfrak{g}$ with an involution $*$ defined on the generators as follows:
$$ (K_j^{\pm1})^* = K_j^{\pm1},\qquad
j=1,2,\ldots,l,
$$
$$
E_j^* = \left\{
\begin{array}{rl}
   K_j F_j, & \; j\neq l_0, \\
  -K_j F_j, & \; j=l_0,
\end{array}
\right.,\qquad F_j^* = \left\{
\begin{array}{rl}
   E_j K_j^{-1}, &\; j \neq l_0, \\
  -E_j K_j^{-1}, &\; j=l_0.
\end{array}
\right.
$$
The Hopf $*$-algebra $(U_q \mathfrak{g},*)$ is a q-analog of the algebra $U
\mathfrak{g}_{\mathbb R} \otimes_{\mathbb R} \mathbb C$ with $\mathfrak{g}_{\mathbb R
}=\operatorname{Lie}(G_{\mathbb R})$. This involution allows us to consider quantum
analogs of noncompact real Lie groups. One can use similar settings to introduce quantum
analogs of compact real groups via the involution $\star$ such that
$$(K_j^{\pm1})^{\star} = K_j^{\pm1},\qquad j=1,2,\ldots,l,
$$
$$
E_j^{\star} = K_j F_j, \qquad F_j^{\star} =  E_j K_j^{-1}, \qquad j=1,..,l.
$$
These two involutions satisfy the relation $*=\theta_q \star=\star \theta_q$ with
$\theta_q$ being a q-analog for the Cartan involution in $U_q \mathfrak{g}$
$$\theta_q(K_j^{\pm1}) = K_j^{\pm1},\qquad j=1,2,\ldots,l,
$$
$$
\theta_q(E_j) = \left\{
\begin{array}{rl}
   E_j,  & \; j\neq l_0, \\
  -E_j, & \; j=l_0,
\end{array}
\right.,\qquad \theta_q(F_j) = \left\{
\begin{array}{rl}
   F_j, &\; j \neq l_0, \\
  -F_j, &\; j=l_0.
\end{array}
\right.
$$

Now we can introduce an antilinear involutive anti-homomorphism $*$ in ${\rm Pol}({\frak
p}^-)_q$ such that ${\rm Pol}({\frak p}^-)_q$ becomes a $(U_q \mathfrak{g},*)$-module
$*$-algebra, i.e.
$$
(\xi f)^*=S(\xi)^*f^*, \qquad f \in {\rm Pol}({\frak p}^-)_q, \xi \in U_q \mathfrak{g}.
$$

\begin{lemma}\label{finite_deg}
For any $\xi \in U_q\mathfrak{g}$ there exist $N \in \mathbb{N}$ such that
$$
\xi \left(\mathbb{C}[\mathfrak{p}^-]_{q,k} \; \mathbb{C}[\mathfrak{p}^+]_{q,-j} \right)
\subset \bigoplus_{|k'-k| \leq N \; \& \; \ |j'-j| \leq N}
\mathbb{C}[\mathfrak{p}^-]_{q,k'} \; \mathbb{C}[\mathfrak{p}^+]_{q,-j'}
$$
for all $k,j \in \mathbb{Z}_+$.
\end{lemma}

{\bf Proof.} A set of all $\xi$ that satisfy the statement is a subalgebra which contains
all $E_i$, $F_i$, $K_i^\pm$, $i=1,2,\ldots,l$. \hfill $\square$

\begin{lemma}\label{z1}
1. Consider the graded $U_q\mathfrak{g}$-module algebra $\mathbb{C}[\mathfrak{p}^-]_q$.
There exists a unique element $z_{\operatorname{low}} \in
\mathbb{C}[\mathfrak{p}^-]_{q,1}$ such that
\begin{equation}\label{z_l0}
F_{l_0} z_{\operatorname{low}} = q_{l_0}^{1/2},\qquad F_j z_{\operatorname{low}} = 0,\; j
\neq l_0.
\end{equation}
2. One has
\begin{equation}\label{z_l0_weight}
K_i^{\pm 1} z_{\operatorname{low}} = q_i^{\pm a_{i l_0}} z_{\operatorname{low}},\quad
i=1,2,\ldots,l,
\end{equation}
\begin{equation*}
E_{l_0} z_{\operatorname{low}} = -q_{l_0}^{1/2} z_{\operatorname{low}}^2,
\end{equation*}
\medskip
$$
E_{l_0} z_{\operatorname{low}}^* = q_{l_0}^{-3/2},\quad F_{l_0} z_{\operatorname{low}}^*
= -q^{5/2} z_{\operatorname{low}}^{*2},\quad K_i^{\pm 1} z_{\operatorname{low}}^* =
q_i^{\mp a_{i l_0}} z_{\operatorname{low}}^*, \; i=1,2,\ldots,l.
$$
\end{lemma}

{\bf Proof.}
 1. It follows from the definitions that the
$U_q\mathfrak{k}$-module $\mathbb{C}[\mathfrak{p}^-]_{q,1}$ is simple. Hence its lowest
weight vector $z_{\operatorname{low}} \neq 0$ is determined up to multiplication by a
non-zero constant. We claim that
\begin{equation*}
F_{l_0}z_{\operatorname{low}} \neq 0.
\end{equation*}
Indeed, otherwise one has $U_q\mathfrak{n}^- z_{\operatorname{low}} \subset \mathbb{C}
z_{\operatorname{low}}$, where $U_q \mathfrak{n} \subset U_q \mathfrak{g}$ is the
subalgebra generated by all $F_j$, $j=1,...,l$. Thus for any $\xi \in U_q\mathfrak{n}^-$
the linear functional $\xi\, z_{\operatorname{low}} \in \mathbb{C}[\mathfrak{p}^-]_{q,1}
$ annihilates the vector $v(\mathfrak{q}^+,0) \in N(\mathfrak{q}^+,0)$. So
$z_{\operatorname{low}}$ annihilates all vectors from $U_q\mathfrak{n}^-\,
v(\mathfrak{q}^+,0)$. Hence $z_{\operatorname{low}}=0$, which contradicts the choice of
$z_{\operatorname{low}}$.

The element $z_{\operatorname{low}}$ can be chosen in a such way that $F_{l_0}
z_{\operatorname{low}} = q_{l_0}^{1/2}$ and it is the unique element with the properties
(\ref{z_l0}).

 2. \eqref{z_l0_weight} is evident since the lowest weight of the  $U_q\mathfrak{k}$-module
$\mathbb{C}[\mathfrak{p}^-]_{q,1}$ is $\alpha_{l_0}$ while the highest weight vector of
the dual $U_q\mathfrak{k}$-module $N(\mathfrak{q}^+,0)_{-1}$ is $F_{l_0}
v(\mathfrak{q}^+,0)$ with weight $-\alpha_{l_0}$.

The relation $E_{l_0}z_{\operatorname{low}}=\mathrm{const}\cdot z_{\operatorname{low}}^2$
follows from the fact that the lowest weight subspace of the $U_q\mathfrak{k}$-module
$\mathbb{C}[\mathfrak{p}^-]_{q,2}$ is generated by $z_{\operatorname{low}}^2$, together
with the fact that $E_{l_0}z_{\operatorname{low}}\in\mathbb{C}[\mathfrak{p}^-]_{q,2}$
belongs to that weight subspace since its weight is  $2\alpha_{l_0}$). The final equality
$\mathrm{const}=-q_{l_0}^{1/2}$ is obtained from the following:
$$
(E_{l_0}F_{l_0}-F_{l_0}E_{l_0})z_{\operatorname{low}}=
\frac{K_{l_0}-K_{l_0}^{-1}}{q_{l_0}-q_{l_0}^{-1}}z_{\operatorname{low}}.
$$
The compatibility of the involutions implies the rest statements of the lemma. \hfill
$\square$


\begin{lemma}\label{z2}
\begin{equation*}
z_{\operatorname{low}}^* z_{\operatorname{low}} = q_{l_0}^2 z_{\operatorname{low}}
z_{\operatorname{low}}^* + 1 - q_{l_0}^2.
\end{equation*}
\end{lemma}

{\bf Proof.} Recall the multiplicative formula for the universal $R$-matrix, keeping in
mind certain reduced expression of $w_0\in W$. Consider the subgroup
$W_{\mathbb{S}}\subset W$ generated by simple reflections $s_j$, $j\ne l_0$. Fix the
reduced expression as follows: $w_0=w_0^\mathbb{S}\cdot\,^\mathbb{S}w_0$, with
$w_0^\mathbb{S} \in W_{\mathbb S}$,
$$
w_0^\mathbb{S}=s_{i_1}s_{i_2}\ldots s_{i_{M'}},\qquad
\,^\mathbb{S}w_0=s_{i_{M'+1}}s_{i_{M'+2}}\ldots s_{i_M}.
$$
The definition of the multiplication $m:\mathrm{Pol}(\mathfrak{p}^-)_q^{\otimes
2}\to\mathrm{Pol}(\mathfrak{p}^-)_q$ implies that (up to flipping tensors)
$m(z_{\operatorname{low}}^*\otimes z_{\operatorname{low}})$ equals
$$
\exp_{q_{M'}^2}((q_{M'}^{-1})E_{\beta_{M'}}\otimes F_{\beta_{M'}})
q^{-t_0}z_{\operatorname{low}}^*\otimes z_{\operatorname{low}}.
$$
One can easily show that
\begin{equation}\label{two}
\beta_{M'}=\alpha_{l_0},\qquad E_{\beta_{M'}}\otimes F_{\beta_{M'}}= c\cdot
E_{l_0}\otimes F_{l_0},\quad c\in\mathbb{C}\backslash\{0\}.
\end{equation}
Now one has
$$
q^{-t_0}(z_{\operatorname{low}}^*\otimes z_{\operatorname{low}})=
q_{l_0}^{-\frac{H_{l_0}\otimes H_{l_0}}2}(z_{\operatorname{low}}^*\otimes
z_{\operatorname{low}}).
$$
Hence,
\begin{equation}\label{c}
z_{\operatorname{low}}^*z_{\operatorname{low}}=
q_{l_0}^2z_{\operatorname{low}}z_{\operatorname{low}}^*+c(1-q_{l_0}^2),\quad c\ne 0.
\end{equation}

It remains to show that $c=1$. For that it suffices to apply $F_{l_0} \in
U_q\mathfrak{g}$ to both sides of the relation
$$
z_{\operatorname{low}}^* z_{\operatorname{low}}^2 = q_{l_0}^4 z_{\operatorname{low}}^2
z_{\operatorname{low}}^* + c(1-q_{l_0}^4) z_{\operatorname{low}},
$$
which leads to a second-degree equation in $c$, whose solutions ($c=0$, $c=1$) are easily
guessed. \hfill $\square$

\bigskip
\subsection{The simplest example}

Let ${\frak g}=\frak{sl}_2$. The algebra $U_q \frak{sl}_2$ is given by its generators $K^{\pm
1},\,E,\,F$ and the relations
$$KK^{-1}=K^{-1}K=1,\quad K^{\pm 1}E=q^{\pm 2}EK^{\pm 1},\quad K^{\pm
1}F=q^{\mp 2}FK^{\pm 1},$$
$$EF-FE=(K-K^{-1})/(q-q^{-1}).$$

 Recall that the comultiplication $\Delta$, the counit $\varepsilon$, and the antipode $S$
are defined on the above generators as follows:
$$\Delta(E)=E \otimes 1\,+\,K \otimes E,\quad \Delta(F)=F \otimes
K^{-1}\,+\,1 \otimes F,\quad \Delta(K^{\pm 1})=K^{\pm 1} \otimes
K^{\pm 1};$$
$$\varepsilon(E)=\varepsilon(F)=0,\quad \varepsilon(K^{\pm 1})=1;$$
$$S(E)=-K^{-1}E,\quad S(F)=-FK,\quad S(K^{\pm})=K^{\mp}.$$

In the notation
$$q=e^{-h/2},\quad K^{\pm 1}=e^{\mp hH/2},\quad E=X^+e^{-hH/4},\quad
F=e^{hH/4}X^-$$ the V. G. Drinfeld formula for the universal R-matrix \cite{DrinfEng}
takes the form
$$
R={\rm exp}_{q^2}((q^{-1}-q)E \otimes F)\cdot {\rm exp}(H \otimes H \cdot
h/4),
$$
with $\exp_t(x)\:=\:\sum \limits_{n=0}^\infty x^n\left(\prod
\limits_{j=1}^n \frac{\textstyle 1-t^j}{\textstyle 1-t}\right)^{-1}$.

The involution $*$ in $U_q \frak{su}(1,1)=(U_q \frak{sl}_2,^*)$ is
defined on the generators $E,\,F,\,K^{\pm 1}$ by
$$E^*=-KF,\quad F^*=-EK^{-1},\quad (K^{\pm 1})^*=K^{\pm 1}.$$

It is well-known that the $*$-algebra $\operatorname{Pol}(\mathbb C )_q$ is
defined by a generator $z$ and the defining relation
$$
z^*z=q^2zz^*+1-q^2.
$$
The $U_q \mathfrak{sl}_2$-module algebra structure is given by
$$F z = q^{1/2},\quad K_i^{\pm 1} z = q_i^{\pm 2} z, \quad E z = -q^{1/2} z^2,$$
$$ E z^* = q^{-3/2},\quad F_{l_0} z^* =
-q^{5/2} z^{*2},\quad K_i^{\pm 1} z^* = q_i^{\pm 2} z^*.
$$

\subsection{\boldmath An embedding
$\operatorname{Pol}(\mathbb{C})_q\hookrightarrow \mathbb{C}[w_0SU_{1,1}]_{q,x}$}
\label{canon_example_sl_2}

The unit disc $\mathbb{D}=\{z \in \mathbb C| |z|<1\}$ is a $SU_{1,1}$-space and in the
category of $SU_{1,1}$-spaces
\begin{equation}\label{naive_embedding}
\mathbb{D}\cong U_1\setminus SU_{1,1},
\end{equation}
where $U_1=\left\{\mathrm{diag}(e^{i\varphi},e^{-i\varphi})|\:
\varphi\in\mathbb{R}/(2\pi\mathbb{Z})\right\}$. The isomorphism
\eqref{naive_embedding} yields an embedding of the $*$-algebra
$\mathrm{Pol}(\mathbb{C})$ into the $*$-algebra of rational
functions on the real affine algebraic group $SU_{1,1}$ such that
\begin{equation}\label{embed_SU11}
z\mapsto t_{11}^{-1}t_{12}.
\end{equation}
Does this embedding have a q-analog? To obtain a positive answer,
one must change the question slightly. Let
$w_0=\begin{pmatrix}0, & -1\\ 1, & 0\end{pmatrix}$. The group
$SU_{1,1}$ is to be substituted by its principal homogeneous space
\begin{equation}\label{su_1_1-orbit}
w_0\,SU_{1,1}=\left\{\left.\left(\begin{array}{cl}t_{11}, & \quad t_{12}\\
t_{21}, & \quad t_{22}\end{array}\right)\in
SL_2\quad\right|\quad\overline{t}_{11}=-t_{22},\:\overline{t}_{12}=
-t_{21}\right\}.
\end{equation}
and the embedding \eqref{embed_SU11} by the embedding
\begin{equation}\label{embed_hat}
z\mapsto t_{12}^{-1}t_{11}.
\end{equation}

\bigskip
We now introduce a quantum analog of the algebra of regular functions on
the $SU_{1,1}$-orbit \eqref{su_1_1-orbit}: equip ${\mathbb
C}[SL_2]_q$ with the involution
\begin{equation}\label{inv}
t_{11}^*=-t_{22},\qquad t_{12}^*=-qt_{21}
\end{equation}
and denote by $\mathbb{C}[w_0SU_{1,1}]_q=(\mathbb{C}[SL_2]_q,*)$ the
defined $*$-algebra. To justify such an involution, note that under the
formal passage to a limit as $q\to 1$ one gets a system of equations
$\bar{t}_{11}=-t_{22}$, $\bar{t}_{12}=-t_{21}$. This allows one to recover
the real space $w_0SU_{1,1}$ from $SL_2$. A more substantial reason is
Proposition \ref{*-hat} below.

Recall that the algebra $\mathbb{C}[SL_2]_q$ of regular functions on the quantum group
$SL_2$ is generated by the matrix elements $t_{ij}$ of the vector representation of
$U_q\mathfrak{sl}_2$, and ${\mathbb C}[SL_2]_q$ is equipped with a
$U_q\mathfrak{sl}_2$-module algebra structure, see \cite{KlSch}.

\bigskip The next statement follows from the definitions.
\begin{proposition}\label{*-hat}
$\mathbb{C}[w_0SU_{1,1}]_q$ is a $U_q\mathfrak{su}_{1,1}$-module
algebra, i.e. one has
\begin{equation}\label{*_*}
(\xi f)^*=(S(\xi))^*f^*,\qquad \xi\in U_q\mathfrak{su}_{1,1},\;
f\in\mathbb{C}[w_0SU_{1,1}]_q.
\end{equation}
\end{proposition}

\medskip

As in the classical case $q=1$, cf. \eqref{embed_hat}, the algebra of
regular functions is too small for our purposes, so we have to extend it.
Recall the useful fact that the element
\begin{equation*}
x=-qt_{12}t_{21}
\end{equation*}
is self-adjoint and quasicommutes with all $t_{ij}$:
$$
t_{11}\,x=q^2\,x\,t_{11},\qquad t_{22}\,x=q^{-2}\,x\,t_{22},\qquad
t_{12}\,x=x\,t_{12}, \qquad t_{21}\,x=x\,t_{21}.
$$
Hence the multiplicative system $x^{\mathbb{Z}_+}$ satisfies the Ore condition. It is
well-known that $\mathbb{C}[SL_2]_q$ is an integral domain, see \cite[p. 266]{Jo}.
Therefore it admits an embedding into its own localization
$\mathbb{C}[w_0SU_{1,1}]_{q,x}$ with respect to the multiplicative system
$x^{\mathbb{Z}_+}$. The self-adjointness of $x$ allows us to extend the involution $*$ on
$\mathbb{C}[w_0SU_{1,1}]_{q,x}$ and to prove that the natural embedding
$\mathbb{C}[w_0SU_{1,1}]_{q,x}\hookrightarrow\mathbb{C}[w_0SU_{1,1}]_{q,x}$ is a
$*$-algebra homomorphism.

Evidently, $t_{12}$, $t_{21}$ are invertible in
$\mathbb{C}[w_0SU_{1,1}]_{q,x}$:
$$t_{12}^{-1}=t_{21}\,x,\qquad t_{21}^{-1}=t_{12}\,x.$$

\begin{proposition}\label{*_local}
\begin{itemize}
\item[1.] The
$U_q\mathfrak{su}_{1,1}$-module algebra structure is uniquely
extendable from $\mathbb{C}[w_0SU_{1,1}]_q$ to
$\mathbb{C}[w_0SU_{1,1}]_{q,x}$.

\item[2.] For any $\xi\in U_q\mathfrak{su}_{1,1}$,
$f\in\mathbb{C}[w_0SU_{1,1}]_{q,x}$ there exists a unique Laurent
polynomial $p_{f,\xi}$ in one variable with coefficients from
$\mathbb{C}[w_0SU_{1,1}]_{q,x}$ such that
\begin{equation}\label{p_f}
\xi(f\cdot x^n)=p_{f,\xi}(q^n)\cdot x^n
\end{equation}
for all $n\in\mathbb{Z}$.
\end{itemize}
\end{proposition}

{\bf Proof.} The uniqueness is obvious in both cases. It is easy to prove that for any
polynomial $\psi$ in one variable
\begin{equation*}
K^{\pm 1}\psi(x)=\psi(x),
\end{equation*}
\begin{equation*}
E\psi(x)=-q^{1/2}t_{11}\frac{\psi(q^{-2}x)-\psi(x)}{q^{-2}x-x}t_{21}=
-q^{1/2}t_{11}t_{21}\frac{\psi(q^{-2}x)-\psi(x)}{q^{-2}x-x},
\end{equation*}
\begin{equation*}
F\psi(x)=-q^{3/2}t_{12}\frac{\psi(q^{-2}x)-\psi(x)}{q^{-2}x-x}t_{22}=
-q^{3/2}t_{12}t_{22}\frac{\psi(x)-\psi(q^2x)}{x-q^2x}.
\end{equation*}
It follows from these observations that for any $\xi\in U_q\mathfrak{su}_{1,1}$,
$f\in\mathbb{C}[w_0SU_{1,1}]_q$ there exists a Laurent polynomial $p_{f,\xi}$ with
coefficients from $\mathbb{C}[w_0SU_{1,1}]_{q,x}$ such that \eqref{p_f} holds for all
$n\in\mathbb{Z}_+$. Any Laurent polynomial is uniquely determined by its values at the
points of the geometric series $q^\mathbb{N}$. This fact implies that \eqref{p_f}
determines a well-defined action of $\xi\in U_q\mathfrak{su}_{1,1}$ on
$\mathbb{C}[w_0SU_{1,1}]_{q,x}$ since
$$p_{f,\xi}(q^n)x=p_{f\cdot x,\xi}(q^{n-1}),\qquad n\in\mathbb{Z}.$$
One can easily prove that the corresponding map
$$
\pi:U_q\mathfrak{su}_{1,1}\to
\operatorname{End}(\mathbb{C}[w_0SU_{1,1}]_{q,x})
$$
is linear and an algebra homomorphism. Furthermore, the involutions \eqref{*_*} are
compatible for obvious reasons. Indeed, the required identities can easily be written in
terms of Laurent polynomials by replacing $f\in\mathbb{C}[w_0SU_{1,1}]_{q,x}$ by the
elements $f_n=f\cdot x^n$, $n\in\mathbb{Z}_+$. It remains to note that any of the Laurent
polynomial identities holds at the points $q^n$ for $n\in\mathbb{N}$ large enough.

Finally, we consider the formula
\begin{equation}\label{Umod_extension}
 \xi(f_1f_2)=\sum_i \xi'_i(f_1)\, \xi''_i(f_2),\qquad
 f_1, f_2 \in \mathbb{C}[w_0SU_{1,1}]_{q,x},\; \xi \in
  U_q\mathfrak{sl}_2,
\end{equation}
with $\Delta \xi=\sum\limits_i \xi'_i \otimes \xi''_i$. As before,
one can prove the existence and uniqueness of a Laurent
polynomial $p_{f,\xi}$ in two variables such that for all $m, n$
$$ \xi(x^m f x^n)\,=\,x^m\, p_{f,\xi}(q^m,q^n)\, x^n.$$
For $m,n \in \mathbb{Z}_+$ large enough, we have
\begin{equation}\label{two_sided_x}
 x^{-m}\,\xi(x^mf_1f_2x^n)\, x^{-n}= \sum_i (x^{-m}\xi'_i(x^mf_1))\,
(\xi''_i(f_2x^n)\,x^{-n})
\end{equation}
since $\mathbb{C}[SL_2]_q$ is a $U_q\mathfrak{sl}_2$-module
algebra. Thus we get the equality of two Laurent polynomials at
all the points $\{(q^m,q^n)\}$ with  $m,n$ large enough. Hence
these polynomials are equal and \eqref{two_sided_x} holds for all
$m,n \in \mathbb{Z}$. Fix $m=n=0$ to get \eqref{Umod_extension}.
\hfill $\square$

\begin{proposition}\label{*_embed}
The map $\mathcal{I}:z\mapsto t_{12}^{-1}t_{11}$ can be uniquely
extended to an embedding of $U_q\mathfrak{su}_{1,1}$-module
algebras
\begin{equation*}
\mathcal{I}:\operatorname{Pol}(\mathbb{C})_q\hookrightarrow
\mathbb{C}[w_0SU_{1,1}]_{q,x}.
\end{equation*}
\end{proposition}

{\bf Proof.} The uniqueness of the $*$-algebra homomorphism is obvious,
while its existence is due to the following identities in
$\mathbb{C}[w_0SU_{1,1}]_{q,x}$:
$$(t_{12}^{-1}t_{11})^*=q^{-1}t_{22}t_{21}^{-1},\quad
(t_{12}^{-1}t_{11})(t_{22}t_{21}^{-1})-
q^{-2}(t_{22}t_{21}^{-1})(t_{12}^{-1}t_{11})=q-q^{-1}.$$

To prove the uniqueness of $\mathcal{I}$, consider a pre-Hilbert space with an
orthonormal basis $\{e_j\}_{j\in\mathbb{Z}_+}$ and a representation $\pi_+$ of
$\mathbb{C}[SL_2]_q$ such that
\begin{align*} \pi_+(t_{11})e_j&=\sqrt{q^{-2(j+1)}-1}\,e_{j+1}, &
\pi_+(t_{12})e_j&=q^{-j}e_j,
\\ \pi_+(t_{21})e_j&=-q^{-(j+1)}e_j, & \pi_+(t_{22})e_j&=
\begin{cases}
-\sqrt{q^{-2j}-1}\,e_{j-1}, & j\ne 0,
\\ 0, & j=0.
\end{cases}
\end{align*}
This is a $*$-representation of $\mathbb{C}[w_0SU_{1,1}]_q$. Maintain the notation
$\pi_+$ for its extension to a $*$-representation of $\mathbb{C}[w_0SU_{1,1}]_{q,x}$.

It follows from the definitions that the representation $\pi_+\circ\mathcal{I}$ of
$\operatorname{Pol}(\mathbb{C})_q$ is equivalent to its Fock representation. Hence
$\pi_+\circ\mathcal{I}$ is faithful. Therefore $\mathcal{I}$ is injective.

Now let us prove that $\mathcal{I}$ is a morphism of $U_q\mathfrak{sl}_2$-modules. Since
$\operatorname{Pol}(\mathbb{C})_q$ and $\mathbb{C}[w_0SU_{1,1}]_{q,x}$ are
$U_q\mathfrak{su}_{1,1}$-module algebras, it is sufficient to prove $\mathcal{I}(\xi
f)=\xi\mathcal{I}(f)$ for $\xi\in U_q\mathfrak{su}_{1,1}$, $f\in\mathbb{C}[z]_q$. For
what remains, we use the embedding of $U_q\mathfrak{sl}_2$-module algebras
$$
\mathbb{C}[t_1,t_2]_q \to\mathbb{C}[SL_2]_q,\qquad t_1\mapsto t_{11},\;t_2 \mapsto t_{12}
$$
(here $\mathbb C[t_1,t_2]_q$ is the well-known algebra with two generators $t_1, t_2$ and
the relation $t_1t_2=qt_2t_1$, see \cite[p.101]{KlSch}) and the embedding of the
$U_q\mathfrak{sl}_2$-module algebra $\mathbb{C}[z]_q$ into $\mathbb{C}[t_1,t_2]_{q,t_2}$
defined as follows:
$$z \mapsto t_2^{-1}t_1. \eqno\square$$

\bigskip

The described embedding of the $U_q\mathfrak{su}_{1,1}$-module algebra
$\operatorname{Pol}(\mathbb{C})_q$ into $\mathbb{C}[w_0SU_{1,1}]_{q,x}$ will be called the canonical embedding.

\bigskip

\begin{remark}\label{SU_2}

Example: the algebra of functions on the quantum group $SU_2$.

Consider the well-known Hopf algebra $U_q\mathfrak{sl}_2$ and the Hopf
$*$-algebra $U_q\mathfrak{su}_2=(U_q\mathfrak{sl}_2,\star)$ together with
the Hopf algebra $\mathbb{C}[SL_2]_q$.

Equipping $\mathbb{C}[SL_2]_q$ with an involution by duality, we obtain the
Hopf $*$-algebra $\mathbb{C}[SU_2]_q$, which is a quantum analog to the
algebra of regular functions on the group $SU_2$. The definitions imply
\begin{equation}\label{ast2}
\left(
\begin{array}{cc}
t_{11}^\star, & t_{12}^\star
\\ t_{21}^\star, & t_{22}^\star
\end{array}
\right)= \left(
\begin{array}{cc}
t_{22}, & -qt_{21}
\\ -q^{-1}t_{12}, & t_{11}
\end{array}
\right).
\end{equation}

Now recall an irreducible $*$-representation of $\mathbb{C}[SU_2]_q$. Let
$\{e_j\}_{j=0}^\infty$ be the standard basis of the Hilbert space $l^2(\mathbb{Z}_+)$.
The following defines an irreducible infinite dimensional $*$-representation $\Pi$ of
$\mathbb{C}[SU_2]_q$ in $l^2(\mathbb{Z}_+)$:
$$ \Pi(t_{11})e_j=\begin{cases}
\sqrt{1-q^{2j}}\,e_{j-1},\;&\; j>0,\\ 0,\;&\;j=0
\end{cases},
 \qquad \Pi(t_{12})e_j=q^{j+1}e_j,$$
 $$ \Pi(t_{22})e_j=\sqrt{1-q^{2(j+1)}}\,e_{j+1}, \qquad \qquad \qquad
  \Pi(t_{21})e_j=q^j e_j.
$$
\end{remark}

\medskip

\subsection{A Fock representation}\label{finite}

In this section we recall the notion of a Fock representation of the algebra
$\mathrm{Pol}(\mathfrak{p}^-)_q$. It was introduced in \cite{SSV4} in the special case
$\mathfrak{p}^-=\mathrm{Mat}_n$. Also it was proved there that the Fock representation is
a unique (up to a unitary equivalence) irreducible faithful representation by bounded
operators in a Hilbert space. We now introduce the Fock representation of the
\hbox{$*$-algebra} $\mathrm{Pol}(\mathfrak{p}^-)_q=
(\mathbb{C}[\mathfrak{p}^-\oplus\mathfrak{p}^+]_q, *)$ for the general case although we
do not discuss many important problems related to it.

Consider the $\mathbb{C}[\mathfrak{p}^- \oplus \mathfrak{p}^+]$-module $\mathcal{H}$ with
one generator $v_0$ and the defining relations
\begin{equation}\label{Fock1}
f v_0 = 0,\qquad f \in \bigoplus\limits_{j=1}^\infty
\mathbb{C}[\mathfrak{p}^+]_{q,-j},
\end{equation}
or, equivalently,
\begin{equation}\label{Fock2}
  f v_0 = 0,\qquad f \in \mathbb{C}[\mathfrak{p}^+]_{q,-1}.
\end{equation}

\begin{lemma}\label{H_as_vect}
The linear map
\begin{equation*}
  \mathbb{C}[\mathfrak{p}^-]_q \rightarrow \mathcal{H},\quad f\mapsto f v_0
\end{equation*}
is an isomorphism of vector spaces.
\end{lemma}

{\bf Proof.} The kernel of the linear map $\mathbb{C}[\mathfrak{p}^- \oplus
\mathfrak{p}^+]_q \rightarrow \mathcal{H}$, $f \mapsto f v_0$ is a left ideal of the
algebra $\mathbb{C}[\mathfrak{p}^- \oplus \mathfrak{p}^+]_q$ generated by
$\bigoplus\limits_{j=1}^\infty \mathbb{C}[\mathfrak{p}^+]_{q,-j}$. This ideal coincides
with $\bigoplus\limits_{i=0}^\infty \bigoplus\limits_{j=1}^\infty
\mathbb{C}[\mathfrak{p}^-]_{q,i} \otimes \mathbb{C}[\mathfrak{p}^+]_{q,-j}$. Hence
$\mathcal{H} \cong \bigoplus\limits_{i=0}^\infty \mathbb{C}[\mathfrak{p}^-]_{q,i} =
\mathbb{C}[\mathfrak{p}^-]_q$. \hfill $\square$

\medskip

The vector space $\mathrm{Pol}(\mathfrak{p}^-)_q$ decomposes as follows:
\begin{equation}\label{expan1}
\mathrm{Pol}(\mathfrak{p}^-)_q = \bigoplus\limits_{i,j=0}^\infty
\mathbb{C}[\mathfrak{p}^-]_{q,i} \cdot \mathbb{C}[\mathfrak{p}^+]_{q,-j},
\end{equation}
and $\mathbb{C}[\mathfrak{p}^-]_{q,0} \cdot \mathbb{C}[\mathfrak{p}^+]_{q,0} = \mathbb{C}
\cdot 1$. That is, for any $f \in \mathrm{Pol}(\mathfrak{p}^-)_q$ there exists a unique
expansion
\begin{equation}\label{expan2}
f = \langle f \rangle \cdot 1 + \sum\limits_{\{ (i,j) \in \mathbb{Z}_+^2|(i,j)\neq
(0,0)\}} f_{i,j},\qquad f_{i,j} \in \mathbb{C}[\mathfrak{p}^-]_{q,i} \cdot
\mathbb{C}[\mathfrak{p}^+]_{q,-j},
\end{equation}
where $\langle \cdot \rangle$ is the linear functional on
$\mathrm{Pol}(\mathfrak{p}^-)_q$ determined by the last formula.

\begin{proposition}\label{inner1}
The sesquilinear form $(f_1,f_2) \stackrel{\operatorname{def}}{=} \langle f_2^* f_1
\rangle$ in $\mathbb{C}[\mathfrak{p}^-]_q$ is non-degenerate, and
\begin{equation}\label{Re}
(f_1, f_2) = \overline{(f_2, f_1)},\qquad f_1,f_2 \in \mathbb{C}[\mathfrak{p}^-]_q.
\end{equation}
\end{proposition}

{\bf Proof.} Firstly, we prove that
\begin{equation}\label{inner_new}
\langle f_2^*f_1\rangle=f_2^*\otimes f_1(R(v(\mathfrak{q}^-,0)\otimes
v(\mathfrak{q}^+,0))),\qquad f_1,f_2\in\mathbb{C}[\mathfrak{p}^-]_q,
\end{equation}
where $R$ is the universal $R$-matrix. Recall that
$$
N(\mathfrak{q}^-,0)\otimes N(\mathfrak{q}^+,0)= \bigoplus\limits_{i,j=0}^\infty
N(\mathfrak{q}^-,0)_i\otimes N(\mathfrak{q}^+,0)_{-j}.
$$
Equip the finite dimensional vector spaces $N(\mathfrak{q}^-,0)_i\otimes
N(\mathfrak{q}^+,0)_{-j}$ with the standard topology, and $N(\mathfrak{q}^-,0)\otimes
N(\mathfrak{q}^+,0)$ with the weakest topology where all projections on
$N(\mathfrak{q}^-,0)_i\otimes N(\mathfrak{q}^+,0)_{-j}$ along to $\bigoplus\limits_{i'\ne
i,j'\ne j}\;N(\mathfrak{q}^-,0)_{i'}\otimes N(\mathfrak{q}^+,0)_{-j'}$ are continuous.
The space of formal series
$$
\sum\limits_{i,j} f_{i,j},\qquad f_{i,j}\in{N(\mathfrak{q}^-,0)_i\otimes
N(\mathfrak{q}^+,0)_{-j}}
$$
is a completion of the topological vector space $N(\mathfrak{q}^-,0)\otimes
N(\mathfrak{q}^+,0)$ and is denoted by $N(\mathfrak{q}^-,0)\widehat{\otimes}
N(\mathfrak{q}^+,0)$. The pairing
$(\mathbb{C}[\mathfrak{p}^+]_q\otimes\mathbb{C}[\mathfrak{p}^-]_q)\times
(N(\mathfrak{q}^-,0)\otimes N(\mathfrak{q}^+,0))\to\mathbb{C}$ extends by continuity to a
pairing $(\mathbb{C}[\mathfrak{p}^+]_q\otimes\mathbb{C}[\mathfrak{p}^-]_q)\times
(N(\mathfrak{q}^-,0)\widehat{\otimes}N(\mathfrak{q}^+,0))\to\mathbb{C}$. Since
$$
\varphi=\sum\limits_i\varphi_i,\quad\varphi_i\in
\mathbb{C}[\mathfrak{p}^-]_{q,i},\qquad\quad\psi=\sum\limits_i\psi_j,\quad
\psi_j\in\mathbb{C}[\mathfrak{p}^+]_{q,-j},
$$
we have $\varphi_0=\varphi(v(\mathfrak{q}^+,0))$, $\psi_0=\psi(v(\mathfrak{q}^-,0))$.
Hence the definition of the multiplication in
$\mathbb{C}[\mathfrak{p}^-\oplus\mathfrak{p}^+]_q$ and $S\otimes S(R)=R$ imply
\eqref{inner_new}.

The action of $U_q\mathfrak{g} \otimes U_q\mathfrak{g}$ in $N(\mathfrak{q}^-,0) \otimes
N(\mathfrak{q}^+,0)$ and the antilinear operator $v_1 \otimes v_2 \mapsto v_2^* \otimes
v_1^*$ in this space are extendable by continuity to $N(\mathfrak{q}^-,0)
\widehat{\otimes} N(\mathfrak{q}^+,0)$.

Since \eqref{inner_new}, the equality \eqref{Re} follows from
$$
R^{21}\, (v(\mathfrak{q}^+,0)\otimes v(\mathfrak{q}^-,0))= \left(R\,
(v(\mathfrak{q}^-,0)\otimes v(\mathfrak{q}^+,0))\right)^{*\otimes *}.
$$
The latter can be checked explicitly.

The non-degeneracy of the sesquilinear form $(\cdot,\cdot)$ in
$\mathbb{C}[\mathfrak{p}^-]_q$ follows from \eqref{inner_new} and the multiplicative
formula for the universal $R$-matrix (see \cite[p. 259]{KlSch}). \hfill $\square$

\medskip

\begin{proposition}\label{preHilbert}
The Hermitian form
$$
(f_1v_0,f_2v_0)=\langle f_2^*f_1\rangle,\qquad f_1,f_2\in
\mathrm{Pol}(\mathfrak{p}^-)_q,
$$
in $\mathcal{H}$ is well defined, positive definite, and one has
$$
(fv',v'')=(v',f^*v''),\qquad v',v''\in\mathcal{H},\;f\in
\mathrm{Pol}(\mathfrak{p}^-)_q.
$$
\end{proposition}

\bigskip

The last proposition allows us to introduce a representation $T_F$ of the $*$-algebra
$\mathrm{Pol}(\mathfrak{p}^-)_q$ in the pre-Hilbert space $\mathcal{H}$. It is called the
Fock representation, and the vector $v_0$ is called the vacuum vector.

\medskip

Consider the special case $\mathfrak{g}=\mathfrak{sl}_2$. It is easy to show that
$$
(z^nv_0,z^nv_0)=((1-q^2y)z^{n-1}v_0,z^{n-1}v_0)=
(1-q^{2n})(z^{n-1}v_0,z^{n-1}v_0)=(q^2;q^2)_n.
$$
Hence the action of $T_F(z)$, $T_F(z^*)$ in the orthonormal basis
$\left\{\sqrt{(q^2;q^2)_n}\;z^n\,v_0\right\}_{n\in\mathbb{Z}_+}$
is described via
\begin{equation*}
T_F(z)e_n=\sqrt{1-q^{2(n+1)}}\, e_{n+1},\qquad T_F(z^*)e_n=
\begin{cases}
\sqrt{1-q^{2n}}\,e_{n-1}, & n\in\mathbb{N},
\\ 0, & n=0.
\end{cases}
\end{equation*}

\medskip

\begin{remark}\label{irrep_Fock}
Obviuosly, the Fock representation is irreducible. Indeed, let
$\mathcal{H}_1\subset\mathcal{H}$ be a non-zero common invariant
subspace of the operators of this representation. Thus
$\mathcal{H}_1$ contains a non-zero vector $v$ such that
$T_F(f)v=0$ for all $f\in\mathbb{C}[\mathfrak{p}^+]_{q,-j}$, $j\ne
0$. This vector is orthogonal to
$\bigoplus\limits_{j=0}^\infty\mathbb{C}[\mathfrak{p}^-]_{q,j}v_0$
, so it belongs to $\mathbb{C}v_0$. Therefore
$v_0\in\mathcal{H}_1$ and $\mathcal{H}_1=\mathcal{H}$.
\end{remark}

\medskip

We now prove the faithfulness of the Fock representation of
$\mathrm{Pol}(\mathfrak{p}^-)_q$. One can prove that
\begin{lemma}\label{ball_I_8.7}
For all $j,k \in \mathbb{Z}_+$ the map
$$
\mathbb{C}[\mathfrak{p}^-]_{q,k} \cdot
\mathbb{C}[\mathfrak{p}^+]_{q,-j} \rightarrow
\mathrm{Hom}(\mathcal{H}_j, \mathcal{H}_k),\qquad f \mapsto
T_F(f)|_{\mathcal{H}_j}
$$
is a bijection.
\end{lemma}

\begin{proposition}\label{ball_I_8.8}
The Fock representation of $\mathrm{Pol}(\mathfrak{p}^-)_q$
is faithful.
\end{proposition}
{\bf Proof.} Equip $\mathbb{Z}_+^2$ with the standard partial
order (i.e. $(k_1,j_1)\le(k_2, j_2)$ if and only if $k_1\le k_2\
\&\ j_1\le j_2$), and the vector space
$\mathrm{Pol}(\mathfrak{p}^-)_q$ with the bigrading
$$
\mathrm{Pol}(\mathfrak{p}^-)_q=
\bigoplus\limits_{k,j=0}^\infty\mathbb{C}[\mathfrak{p}^-]_{q,k}\;
\mathbb{C}[\mathfrak{p}^+]_{q,-j}.
$$
Suppose that there exists $f\ne 0$ such that $T_F(f)=0$.
Use the expansion
$$
f=\sum\limits_{k,j=0}^\infty f_{k,j},\qquad f_{k,j}\in
\mathbb{C}[\mathfrak{p}^-]_{q,k}\mathbb{C}[\mathfrak{p}^+]_{q,-j},
$$
to choose a non-zero component $f_{k_0,j_0}$ in the expansion with a
minimal $(k_0,j_0)$. $T_F(f)=0$ implies that
$T_F(f_{k_0,j_0})|_{\mathcal{H}_{j_0}}=0$ since
$T_F(f_{k_0,j_0})|_{\mathcal{H}_{j_0}}$ is a product of the linear operator
$T_F(f)|_{\mathcal{H}_{j_0}}$ and the projection onto $\mathcal{H}_{k_0}$
along to $\bigoplus\limits_{k\ne k_0}\mathcal{H}_k$. It remains to note
that
$$
f_{k_0,j_0}\ne 0,\qquad
f_{k_0,j_0}\in\mathbb{C}[\mathfrak{p}^-]_{q,k_0}
\mathbb{C}[\mathfrak{p}^+]_{q,-j_0},\qquad
T_F(f_{k_0,j_0})|_{\mathcal{H}_{j_0}}=0.
$$
This is a contradiction with Lemma \ref{ball_I_8.7}.
\hfill $\square$

\section{The canonical embedding}\label{canonical_embed}

\subsection{Introduction}\label{canon_example}

In what follows, all the pairs $(\mathfrak{g},\mathfrak{k})$ are assumed to be Hermitian
symmetric. A duality argument allows one to define a Hopf $*$-algebra
$((U_q\mathfrak{g})^*,\sharp)$:
\begin{equation}\label{g_0}
f^\sharp(\xi)=\overline{f((S(\xi))^*)},\qquad
f\in(U_q\mathfrak{g})^*,\;\xi\in U_q\mathfrak{g}.
\end{equation}

It follows from Burnside's theorem that the homomorphism of
algebras $U_q\mathfrak{g}\to\operatorname{End}L(\Lambda)$ is onto
for all $\Lambda\in P_+$, hence
$(\operatorname{End}L(\Lambda))^*\hookrightarrow
(U_q\mathfrak{g})^*$. Since
$$
\sharp:(\operatorname{End}L(\Lambda))^*\to
(\operatorname{End}L(-w_0\Lambda))^*,
$$
the involution $\sharp$ admits a restriction to the subalgebra
$$
\mathbb{C}[G]_q=\sum\limits_{\Lambda\in
P_+}(\operatorname{End}L(\Lambda))^*\subset(U_q\mathfrak{g})^*.
$$

Consider the Hopf $*$-algebra $\mathbb{C}[G_{\mathbb R}]_q=(\mathbb{C}[G]_q,\sharp)$. It
can be viewed as a $q$-analog of the algebra of regular functions on the real affine
algebraic group $G_{\mathbb R}=\{g\in G\,|\,f^\sharp(g)=\overline{f(g)},\quad
f\in\mathbb{C}[G]\}$. It is a well known fact that $G_{\mathbb R}$ acts by automorphisms
of a bounded symmetric domain $\mathbb{D}$ which admits a standard embedding into
$\mathfrak{p}^-$, see \cite{Kor}.

The presence of a distinguished point $0\in\mathbb{D}$ allows one to associate to every
function $f(z)$ on $\mathbb{D}$ a function $f(g^{-1}\cdot 0)$ on $G_{\mathbb R}$. We are
interested in quantum analogs of function algebras on $\mathbb{D}$ and their images under
the above embedding into algebras of functions on $G_{\mathbb R}$.

In the special case $\mathfrak{g}=\mathfrak{sl}_2$, $G_{\mathbb R}$ is the real affine
algebraic group $SU_{1,1}$ ; it is easy to verify that
\begin{equation}\label{inv_su11}
t_{11}^\sharp=t_{22},\qquad t_{12}^\sharp = q t_{21},
\end{equation}
with $\{t_{ij}\}_{i,j=1,2}$ being the standard generators of
$\mathbb{C}[SL_2]_q$. In fact, consider the automorphism $\theta_q^*$ of
the Hopf algebra $\mathbb{C}[SL_2]_q$ which is dual to the automorphism
$\theta_q:U_q\mathfrak{sl}_2\to U_q\mathfrak{sl}_2$,
$$
\theta_q:K^{\pm 1}\mapsto K^{\pm 1},\qquad\theta_q:E\mapsto
-E,\qquad\theta_q:F\mapsto -F.
$$
Its action on the generators $t_{ij}$ is given by
$$
\begin{pmatrix}\theta_q^*(t_{11}), & \; \theta_q^*(t_{12})\\
\theta_q^*(t_{21}), &\;\theta_q^*(t_{22})\end{pmatrix}=
\begin{pmatrix}t_{11}, & -t_{12}\\ -t_{21}, & t_{22}\end{pmatrix}.
$$
What remains is to use the relation $*=\theta\star$, together with the explicit formulas
\begin{equation*} \left(
\begin{array}{cc}
t_{11}^\star, & t_{12}^\star
\\ t_{21}^\star, & t_{22}^\star
\end{array}
\right)= \left(
\begin{array}{cc}
t_{22}, & -qt_{21}
\\ -q^{-1}t_{12}, & t_{11}
\end{array}
\right).
\end{equation*}

Recall that the canonical embedding in the case
$\mathfrak{g}=\mathfrak{sl}_2$ (that is, $\mathbb{D}=\{z\in\mathbb{C}\,|\,
|z|<1\}$) has been described in details in section
\ref{canon_example_sl_2}.

\medskip

\subsection{\boldmath A quantum analog of the Weyl group element
$w_0$}\label{w_0}

The construction of the canonical embedding becomes considerably more intricate
under the passage from the special case $\mathfrak{g}=\mathfrak{sl}_2$ to the
general case, and the passage requires the notion of the quantum Weyl group
described below.

\begin{proposition}\label{s_bar}
There exists a unique linear functional $\overline{s}$ on
$\mathbb{C}[SL_2]_q$ such that
\begin{align*}
\overline{s}(t_{22}f)&=\overline{s}(ft_{11})=0,&
\\ \overline{s}(t_{12}f)&=\overline{s}(ft_{12})=q\overline{s}(f),&
\\ \overline{s}(t_{21}f)&=\overline{s}(ft_{21})=-\overline{s}(f)&
\end{align*}
for all $f\in\mathbb{C}[SL_2]_q$ with the property $f(1)=1$.
\end{proposition}

{\bf Proof.} Let $\{e_j\}$, $j\in\mathbb{Z}_+$, be the standard basis of the Hilbert
space $l^2(\mathbb{Z}_+)$ and $\Pi$ the $*$-representation of the algebra
$\mathbb{C}[SU_2]_q=(\mathbb{C}[SL_2]_q,\star)$ in this space defined as follows:
$$ \Pi(t_{11})e_j=\begin{cases} \sqrt{1-q^{2j}}\,e_{j-1},\;&\;
j>0,\\ 0,\;&\;j=0
\end{cases},
 \qquad \Pi(t_{12})e_j=q^{j+1}e_j,$$
 $$ \Pi(t_{22})e_j=\sqrt{1-q^{2(j+1)}}\,e_{j+1}, \qquad \qquad \qquad
  \Pi(t_{21})e_j=q^j e_j.
$$

The linear functional $\overline{s}(f)=(\Pi(f)e_0,e_0)$ possesses all the required
properties. The uniqueness of $\overline{s}$ follows from the fact that
$\mathbb{C}[SL_2]_q$ is a linear span of the monomials
$$
\{t_{22}^dt_{21}^ct_{12}^bt_{11}^a\;|\;a,b,c,d\in\mathbb{Z}_+\}.\eqno\square
$$

\medskip

\begin{remark}
It is easy to prove that every element of $\mathbb{C}[SL_2]_q$ admits a
unique decomposition
$$
f\,=\,\sum_{j=1}^\infty
t_{22}^j\,f_j(t_{12},t_{21})\,+\,f_0(t_{12},t_{21})\,+\,\sum_{j=1}^\infty
f_{-j}(t_{12},t_{21})\,t_{11}^j,
$$
with each $f_j$ being a polynomial in two variables. This observation provides another
proof of existence for $\overline{s}$:
$$\overline{s}(f)=f_0(q,-1),\qquad f\in\mathbb{C}[SL_2]_q.$$

\bigskip

$\overline{s}\in\mathbb{C}[SL_2]_q^*$ is a $q$-analog of the element
$\left(\begin{smallmatrix}0 & 1\\ -1 & 0\end{smallmatrix}\right)$ of $SL_2$
since
\begin{equation}\label{s_bar_eq}
\begin{pmatrix}\overline{s}(t_{11}) & \overline{s}(t_{12})\\
\overline{s}(t_{21}) & \overline{s}(t_{22})\end{pmatrix}=\begin{pmatrix}0
& q\\ -1 & 0\end{pmatrix}.
\end{equation}
\end{remark}

We equip the vector space $U_q\mathfrak{sl}_2$ with the weakest among those
topologies in which all the irreducible finite dimensional weight
representations
$\pi_\lambda:U_q\mathfrak{sl}_2\to\operatorname{End}(L(\lambda))$,
$\lambda\in\mathbb{Z}_+\overline{\omega}$, are continuous. As these
representations separate elements of $U_q\mathfrak{sl}_2$, the above
topology is Hausdorff and the completion of $U_q\mathfrak{sl}_2$ is
canonically isomorphic to the direct product of algebras
$\mathop{\times}\limits_j\operatorname{End}L(j\overline{\omega})$, that is,
to the algebra $\mathbb{C}[SL_2]_q^*$ with the weak-$*$ topology. The
following claim is a direct consequence of the definitions.

\begin{lemma}\label{extension_pi}
Every finite dimensional weight representation of the algebra
$U_q\mathfrak{sl}_2$ in a vector space $L$ extends by continuity to a
representation of $\mathbb{C}[SL_2]_q^*$, thus equipping $L$ with the
structure of a $\mathbb{C}[SL_2]_q^*$-module.
\end{lemma}

\bigskip

Recall that, given a weight $U_q\mathfrak{sl}_2$-module $V$, one has well
defined linear operators $H$ and $X^{\pm}$ in $V$ such that
\begin{equation}\label{changesl_2}
K^{\pm1}=q^{\pm H},\quad X^+=Eq^{-\frac{1}{2}H},\quad
X^-=q^{\frac{1}{2}H}F.
\end{equation}

Let $j\in \mathbb{Z}_+$. We follow the conventions of the quantum theory of
angular momentum in introducing the standard basis
$\{e_m^j\}_{m=j,j-1,\ldots,-j}$ of the vector space
$L(2j\overline{\omega})$ by setting
\begin{equation*}
\begin{gathered}
He_m^j=2me_m^j,\quad X^+e_m^j=
\begin{cases}
([j-m]_q[j+m+1]_q)^{1/2}e_{m+1}^j, & m\ne j
\\ 0, & m=j
\end{cases}
\\ X^-e_m^j=
\begin{cases}
([j-m+1]_q[j+m]_q)^{1/2}e_{m-1}^j, & m\ne-j
\\ 0, & m=-j
\end{cases},
\end{gathered}
\end{equation*}
with $[n]_q=\dfrac{q^n-q^{-n}}{q-q^{-1}}$.

\begin{lemma}(\cite[p. 8]{VakSoib88})\label{simple_refl}
For all $j,m$
$$\bar{s}e_m^j=(-1)^{j-m}q^{j^2+j}q^{-(m^2+m)}e_{-m}^j,$$
hence
\begin{equation*}
\overline{s}K\overline{s}^{-1}=K^{-1}.
\end{equation*}
\end{lemma}

\bigskip

Let us turn from the special case $G=SL_2$ to the general case. Recall that
the group $G$ in our situation is assumed to be simply connected. Equip the
vector space $\mathbb{C}[G]_q^*\cong \mathop{\times}\limits_{\lambda\in
P_+}\operatorname{End}(L(\lambda))$ with the direct product topology, i.e.,
the weak-$*$topology. Consider the standard embeddings of Hopf $*$-algebras
$$\varphi_i:\,U_{q_i}\mathfrak{sl}_2\hookrightarrow U_q\mathfrak{g},$$
$$
\varphi_i:K^{\pm 1}\mapsto K_i^{\pm 1},\quad\varphi_i:E\mapsto
E_i,\quad\varphi_i:F\mapsto F_i,\qquad i=1,2,\ldots,l.
$$

\begin{proposition}\label{embed_Sl2}
The embeddings of algebras $\varphi_i:U_{q_i}\mathfrak{sl}_2\hookrightarrow
U_q\mathfrak{g}$ extend by continuity to embeddings of algebras
$\overline{\varphi}_i:\mathbb{C}[SL_2]_{q_i}^*\hookrightarrow
\mathbb{C}[G]_q^*$.
\end{proposition}

{\bf Proof.} The possibility to extend follows from Lemma
\ref{extension_pi}, and the injectivity is due to
$$
\dim\operatorname{Hom}_{\,U_{q_i}\mathfrak{sl}_2}
(L(k\overline{\omega}),L(k\overline{\omega}_i))>0,\qquad k \in\mathbb{Z}_+,
$$
with the $U_q\mathfrak{g}$-module $L(k\overline{\omega}_i)$ being treated
here as a $U_{q_i}\mathfrak{sl}_2$-module. \hfill $\square$

\bigskip

Let $\overline{s}_i=\overline{\varphi}_i(\overline{s})$, $i=1,2,\ldots,l$.
It is well known \cite[p. 157]{LS}, \cite[p. 427]{KR} that the elements
$\overline{s}_1,\overline{s}_2,\ldots,\overline{s}_l$ of the algebra
$\mathbb{C}[G]_q^*$ satisfy the braid relations
\begin{equation}\label{braid_rel}
\overline{s}_i\overline{s}_j\overline{s}_i\cdots=
\overline{s}_j\overline{s}_i\overline{s}_j\cdots,
\end{equation}
where the number of factors in either side of \eqref{braid_rel} is
two if $a_{ij}a_{ji}=0$, three if $a_{ij}a_{ji}=1$, four if
$a_{ij}a_{ji}=2$, and six if $a_{ij}a_{ji}=3$.\footnote{The last
possibility is beyond the Hermitian-symmetric case we consider here.} Let
$w=s_{i_1}s_{i_2}\cdots s_{i_M}$ be a reduced expression of an element
$w\in W$ and set
\begin{equation*}
\overline{w}=\overline{s}_{i_1}\overline{s}_{i_2}\cdots\overline{s}_{i_M}.
\end{equation*}
This element of $\mathbb{C}[G]_q^*$ does not depend on the choice of
reduced expression, as is clear from \eqref{braid_rel}, together with a
well known property of Coxeter groups \cite[p. 19]{Bou4-6}.

We get a quantum analog $\overline{w}$ of $w\in W$ and, in particular, a
quantum analog of the longest element $w_0$ of the Weyl group $W$.

An application of Lemma \ref{simple_refl} allows one to prove

\begin{proposition}(\cite[p. 138]{KorSoib})\label{s_general}
The elements $\overline{s}_i\in\mathbb{C}[G]_q^*$ are invertible, and
$$
\overline{s}_iK_j\overline{s}_i^{-1}=K_jK_i^{-a_{ij}},\qquad
i,j=1,2,\ldots,l.
$$
\end{proposition}

\begin{corollary}\label{refl_weight}
Let $V=\bigoplus\limits_\lambda V_\lambda$ be a finite dimensional
$U_q\mathfrak{g}$-module which may be decomposed as a direct sum of its weight
subspaces. Then
\begin{equation}\label{w_action}
\overline{w}V_\lambda=V_{w\lambda},\qquad w\in W.
\end{equation}
\end{corollary}

\medskip

We claim that the relation
\begin{equation}\label{w0_tilda}
\widetilde{w}_0=\overline{w}_0^{\,-1}\cdot
q^{-\frac12\sum\limits_{i,j=1}^lc_{ij}d_id_jH_iH_j}
\end{equation}
determines an element $\widetilde{w}_0$ of $\mathbb{C}[G]_q^*$.
Firstly, by Proposition \ref{s_general}, the $\overline{s}_i$ are
invertible, hence so is their product $\overline{w}_0$. Secondly,
$\mathbb{C}[G]_q^*$ is a direct product of the finite dimensional
algebras $\operatorname{End}L(\lambda)$ and all $L(\lambda)$ are
$U_q\mathfrak{g}$-weight modules. Hence the second factor in
\eqref{w0_tilda}, which is normally called the Cartan correcting
multiplier, determines an element of $\mathbb{C}[G]_q^*$.

In the case $\mathfrak{g}=\mathfrak{sl}_2$ one has
$\widetilde{w}_0=\bar{s}^{\,-1}\cdot q^{-\frac{H^2}4}$. It can be verified
\cite[p. 425]{KR} that
$$
\widetilde{w}_0^{-1}X^+\widetilde{w}_0=-q^{-1}\,X^-,\qquad
\widetilde{w}_0^{-1}X^-\widetilde{w}_0=-q\,X^+,
$$
\begin{equation*}
\Delta\widetilde{w}_0=(\widetilde{w}_0\otimes\widetilde{w}_0)R,
\end{equation*}
with $\Delta:\mathbb{C}[SL_2]_q^*\to(\mathbb{C}[SL_2]_q^{\otimes 2})^*$
and $R\in(\mathbb{C}[SL_2]_q^{\otimes 2})^*$ corresponding to the universal
$R$-matrix of the Hopf algebra $U_q\mathfrak{sl}_2$, see \cite[p. 425]{KR}.

In a similar way, we introduce in the general case an element $R$ of
$(\mathbb{C}[G]_q^{\otimes 2})^*$ corresponding to the universal $R$-matrix
of the Hopf algebra $U_q\mathfrak{g}$, cf. \cite[p. 289]{Jo}.

\begin{proposition}(\cite[p. 252]{LS})\label{w0_tilda_Prop}
In $(\mathbb{C}[G]_q^{\otimes 2})^*$ the following relation is valid:
\begin{equation}\label{w0_R_new}
\Delta\widetilde{w}_0=(\widetilde{w}_0\otimes\widetilde{w}_0)R
\end{equation}
with $\Delta:\mathbb{C}[G]_q^*\to(\mathbb{C}[G]_q^{\otimes 2})^*$.
\end{proposition}

\begin{remark}\label{to_LS}
In \cite{LS,KR} $\widetilde{w}_0^{-1}$ is used instead of
$\widetilde{w}_0$.
\end{remark}

Consider the elements $X_i^\pm$ of $\mathbb{C}[G]_q^*$ given by
\begin{equation*}
X_i^+=E_iq_i^{-\frac{H_i}{2}},\qquad X_i^-=q_i^{\frac{H_i}{2}}F_i,
\end{equation*}
cf. \eqref{changesl_2}.

The following result is due to Joseph \cite{Jo-R}.
\begin{proposition}\label{from_Jo-R}
In $\mathbb{C}[G]_q^*$ one has
\begin{equation*}
\widetilde{w}_0^{-1}\,X_i^{\pm}\,\widetilde{w}_0\,=\,-q_i^{\mp 1}
X_{i'}^{\mp},\qquad i=1,2,\ldots,l,
\end{equation*}
with $i'$ being determined by $\alpha_{i'}=-w_0\alpha_i$.
\end{proposition}

\begin{remark}\label{to_Joseph}
Our choice of the Cartan correcting multiplier in \eqref{w0_tilda} is
conventional, see, e.g., \cite[p. 252]{LS}, and differs a little from that
used by Joseph, see \cite[p. 422]{Jo-R}.
\end{remark}

\bigskip

\subsection{Fundamental representations and special
bases.}\label{true_bases}

The study of algebras of functions on compact quantum groups extensively utilizes the
matrix elements of representations $\pi_\Lambda$ of $U_q\mathfrak{g}$ corresponding to
the simple finite dimensional $U_q\mathfrak{g}$-weight modules $L(\Lambda)$, $\Lambda\in
P_+$. In this approach, some unpleasantness is related to an ambiguity in choosing an
orthonormal basis of weight vectors in $L(\Lambda)$.

Even in the special case of a fundamental weight $\Lambda$ and $\lambda,\mu\in W\Lambda$,
$\lambda\ne\mu$, with one dimensional weight subspaces, the matrix element
$c^{\Lambda}_{\lambda,\mu}$ is determined only up to a scalar multiplier of modulus one.
We are going to remove this ambiguity by removing the arbitrariness in the choice of an
orthonormal basis of weight vectors.

\bigskip

Let $\mathbb{S}=\{1,2,\ldots,l\}\setminus\{l_0\}$. $W_{\mathbb S}$ is the stabilizer of
the fundamental weight $\overline{\omega}=\overline{\omega}_{l_0}$ \cite[p. 22]{Hum2},
that is $W_{\mathbb{S}}=\{\sigma \in W\,|\,\sigma \overline{\omega}=\overline{\omega}\}$.
Hence the set $W\overline{\omega}$ is in a natural one-to-one correspondence with the set
$^\mathbb{S}W$ of shortest elements in the cosets $\sigma W_\mathbb{S}$, $\sigma \in W$,
see \cite[p. 19]{Hum2}.

Recall that the representation $\pi_{\overline{\omega}}$ of the $*$-algebra
$(U_q\mathfrak{g},\star)$ in $L(\overline{\omega})$ is unitarizable and the Hermitian
form $(\cdot,\cdot)$ in question is uniquely determined by the normalization condition
$(v(\overline{\omega}),v(\overline{\omega}))=1$.

An orthonormal basis formed by weight vectors of the $U_q\mathfrak{g}$-module
$L(\overline{\omega})$ will be called special if the weight vectors with weights in
$W\overline{\omega}$ are given by
$$
v_{\sigma \overline{\omega}}^{\overline{\omega}}=
\frac{\overline{\sigma}^{-1}v(\overline{\omega})}
{\left\|\overline{\sigma}^{-1}v(\overline{\omega})\right\|},\qquad
\sigma \in\,^\mathbb{S}W.
$$
Here $\overline{\sigma}\in\mathbb{C}[G]_q^*$ is the quantum analog of $\sigma$
introduced in the previous section (its invertibility and \eqref{w_action}
are used).

We now give an equivalent definition of a special basis, cf. \cite[p. 153]{LakshResh}.
Let $\sigma=s_{i_L}s_{i_{L-1}}\cdots s_{i_1}$ be a reduced expression for
$\sigma \in\,^\mathbb{S}W$. Set $u_j=s_{i_j}s_{i_{j-1}}\cdots s_{i_1}$,
$$
\lambda_j=
\begin{cases}
u_j\overline{\omega}, & j=1,2,\ldots,L,
\\ \overline{\omega}, & j=0,
\end{cases}\qquad \qquad
m_j=\frac{2(\lambda_j,\lambda_{j-1})}{(\lambda_j,\lambda_j)}.
$$
One has $m_j\in\mathbb{N}$ \cite[p. 93]{Bou4-6}. The validity of the following claim is a
consequence of Lemma \ref{simple_refl}.

\begin{proposition}\label{standard_def}
An orthonormal basis of weight vectors is special if and only if
$$
v_{\sigma \overline{\omega}}^{\overline{\omega}}=
\frac{F_{i_L}^{m_L}F_{i_{L-1}}^{m_{L-1}}\cdots
F_{i_1}^{m_1}v(\overline{\omega})}
{\left\|F_{i_L}^{m_L}F_{i_{L-1}}^{m_{L-1}}\cdots
F_{i_1}^{m_1}v(\overline{\omega})\right\|}
$$
for all $\sigma \in\,^\mathbb{S}W$.
\end{proposition}

\medskip
We introduce the notation
$$
c_{\lambda,i;\mu,j}^\Lambda(\xi)\stackrel{\operatorname{def}}{=}\left(\xi
v_{\mu,j}^\Lambda, v_{\lambda,i}^\Lambda\right)
$$
for matrix elements of the representations $\pi_\Lambda$ corresponding to the
$U_q\mathfrak{g}$-modules $L(\Lambda)$, with $v_{\lambda,i}$ being an orthonormal basis
of weight vectors in any weight subspace $L(\Lambda)_\lambda$. We write
$c_{\lambda,\mu}^\Lambda$ instead of $c_{\lambda,i;\mu,j}^\Lambda$ if this does not lead
to ambiguities.

\begin{proposition}\label{U_standard}
One has the following relations for the matrix elements of the fundamental
representations in the special bases:
\begin{eqnarray}
S(c_{\sigma'\overline{\omega},\sigma''\overline{\omega}}^{\overline{\omega}}) &=&
(-1)^{l(\sigma')-l(\sigma'')}q^{(\mu-\lambda,\rho)}
c_{-\sigma''\overline{\omega},-\sigma'\overline{\omega}}^{-w_0\overline{\omega}}\;,
\label{antipode_standard}
\\ (c_{\sigma'\overline{\omega},\sigma''\overline{\omega}}^{\overline{\omega}})^\star
&=& (-1)^{l(\sigma')-l(\sigma'')}q^{(\lambda-\mu,\rho)}
c_{-\sigma'\overline{\omega},-\sigma''\overline{\omega}}^{-w_0\overline{\omega}}\;,
\label{aster_standard}
\end{eqnarray}
with $\sigma',\sigma''\in\,^\mathbb{S}W$, and $\rho$ being half the sum of positive
roots.
\end{proposition}

{\bf Proof.} Let $\pi_{\overline{\omega}}$, $\pi_{-w_0\overline{\omega}}$ be the
fundamental representations of $U_q\mathfrak{g}$ corresponding to the
$U_q\mathfrak{g}$-modules $L(\overline{\omega})$ and $L(-w_0\overline{\omega})$,
respectively. Equip $U_q\mathfrak{g}$ with the {\bf antilinear} involution $\tau$
\begin{equation*}
\tau(K_j^{\pm1})=S(K_j^{\pm1}),\quad\tau(E_j)=q_jS(E_j),\quad
\tau(F_j)=q_j^{-1}S(F_j),
\end{equation*}
with $j=1,2,\ldots,l$. To any $\Lambda\in P_+$ we assign the
$U_q\mathfrak{g}$-module $L(\Lambda)$ and the $*$-representation
$\pi_\Lambda$ in the finite dimensional Hilbert space
$L(\Lambda)$. Let $\pi_\Lambda^\tau$ be the representation of
$U_q\mathfrak{g}$ in $L(\Lambda)$ defined as follows:
$$
(\pi_\Lambda^\tau(\xi)v_1,v_2)=(v_1,\pi_\Lambda(\tau(\xi))v_2),\qquad \xi\in
U_q\mathfrak{g},\;v_1,v_2\in L(\Lambda).$$ Since $\star\tau=\tau\star$,
$$
(\pi_\Lambda^\tau(\xi)v_1,v_2)=(v_1,\pi_\Lambda(\tau(\xi))v_2)=
(\pi_\Lambda((\tau(\xi))^\star)v_1,v_2)=
$$
$$
=(\pi_\Lambda(\tau(\xi^\star))v_1,v_2)=(v_1,\pi_\Lambda^\tau(\xi^\star)v_2).
$$
Hence $\pi_\Lambda^\tau$ is a $*$-representation with highest weight $-w_0\Lambda$. So we
have shown that there exists a unique isometric linear map
$$
U:L(\overline{\omega})\to L(-w_0\overline{\omega}),\qquad
U:v_{\overline{\omega}}^{\overline{\omega}}\mapsto
v_{-\overline{\omega}}^{-w_0\overline{\omega}},
$$
and this establishes the equivalence of $\pi_{\overline{\omega}}^\tau$ and
$\pi_{-w_0\overline{\omega}}$. Proposition \ref{standard_def} and the definition of
$\tau$ imply that
$$
Uv_{\sigma\overline{\omega}}^{\overline{\omega}}=
(-1)^{l(\sigma)}v_{-\sigma\overline{\omega}}^{-w_0\overline{\omega}},\qquad \sigma\in
^\mathbb{S}W.
$$
What remains is to apply the next proposition which, inessentially, differs from the
Soibelman result, see \cite[p. 100]{KorSoib}.

\begin{proposition}\label{invol_c}
Let $U:L(\Lambda)\to L(-w_0\Lambda)$ be a unique (up to multiplying by a
constant) map such that
$\pi_\Lambda^\tau(\xi)=U^{-1}\,\pi_{-w_0\Lambda}(\xi)\,U$ for all $\xi\in
U_q\mathfrak{g}$. In the $*$-algebra $\mathbb{C}[K]_q$ one has
\begin{equation*}
S\left(c_{\lambda,i;\mu,j}^\Lambda\right)=q^{(\mu-\lambda,\rho)}\cdot
\widetilde{c}_{-\mu,j;-\lambda,i}^{\;-w_0\Lambda},
\end{equation*}
\begin{equation*}
\left(c_{\lambda,i;\mu,j}^\Lambda\right)^\star=q^{(\lambda-\mu,\rho)}\cdot
\widetilde{c}_{-\lambda,i;-\mu,j}^{\;-w_0\Lambda},
\end{equation*}
with $\left\{c_{\lambda,i;\mu,j}^\Lambda\right\}$,
$\left\{\widetilde{c}_{\lambda,i;\mu,j}^{\;-w_0\Lambda}\right\}$
being the matrix elements of the representations $\pi_\Lambda$,
$\pi_{-w_0\Lambda}$ in the bases
$\left\{v_{\mu,j}^\Lambda\right\}$,
$\left\{U\,v_{-\mu,j}^\Lambda\right\}$.
\end{proposition} \hfill $\square$

\bigskip

In closing, consider an important special case of the fundamental weight
$\overline{\omega}$ corresponding to a simple root whose coefficient in the decomposition
of the highest root is $1$. Even more, we restrict ourselves to the simple-laced Lie
algebras (the series $A$, $D$, $E$). Under these assumptions $\overline{\omega}$ is a
microweight, i.e., all the weights of the $U_q\mathfrak{g}$-module $L(\overline{\omega})$
belong to its $W$-orbit \cite[p. 164]{Bou4-6}. Now \eqref{antipode_standard},
\eqref{aster_standard} take the form
\begin{eqnarray*}
S(c_{\lambda,\mu}^{\overline{\omega}}) &=& (-q)^{(\mu-\lambda,\rho)}\,
c_{-\mu,-\lambda}^{-w_0\overline{\omega}},\qquad\lambda,\mu\in W\overline{\omega},
\\ (c_{\lambda,\mu}^{\overline{\omega}})^\star &=&
(-q)^{(\lambda-\mu,\rho)}\,c_{-\lambda,-\mu}^{-w_0\overline{\omega}},\qquad
\lambda,\mu\in W\overline{\omega},
\end{eqnarray*}
because  $l(s_i\sigma)=l(\sigma)+1$, $s_i\sigma\in\,^\mathbb{S}W$, implies that
$$
(s_i\,\sigma\overline{\omega}-\sigma\overline{\omega},\rho)=
-\left(\alpha_i,\sum_{j=1}^l\overline{\omega}_j\right)=-1.
$$

\bigskip

\subsection{\boldmath The canonical embedding
$\mathbb{C}[\mathfrak{p}^-]_q\hookrightarrow
\mathbb{C}[X^-_{\mathbb{S}}]_{q,t}$.}\label{w0_G_x}

We maintain our assumptions that
$\mathbb{S}=\{1,2,\ldots,l\}\setminus\{l_0\}$ and that the simple root
$\alpha_{l_0}$ has coefficient $1$ in the decomposition of the highest
root. Consider the fundamental weight
$$\overline{\omega}_{l_1}=-w_0\overline{\omega}_{l_0}$$
together with the associated fundamental representation $\pi_{l_1}$ of
$U_q\mathfrak{g}$ in the vector space $L(\overline{\omega}_{l_1})$.

Pick a special basis in $L(\overline{\omega}_{l_1})$ and maintain the notation for matrix
elements of representations of $U_q\mathfrak{g}$. We consider the following distinguished
matrix element of $\pi_{l_1}$:
\begin{equation*}
t= c_{\overline{\omega}_{l_1},-\overline{\omega}_{l_0}}^{\overline{\omega}_{l_1}}.
\end{equation*}

One can show

\begin{proposition}\label{t_tstar}
1. $tt^\star=t^\star t$.

2. Let $\Lambda\in P_+$ and $\lambda,\mu$ be weights of $L(\Lambda)$. The
matrix elements $c_{\lambda,j;\mu,k}^{L(\Lambda)}$ of the algebra
$\mathbb{C}[G]_q$ quasi-commute with $t$, $t^\star$. Specifically,
$$
t\,c_{\lambda,j;\mu,k}^{L(\Lambda)}=\mathrm{const}_1\;
c_{\lambda,j;\mu,k}^{L(\Lambda)}\,t,\qquad t^\star\,
c_{\lambda,j;\mu,k}^{L(\Lambda)}=\mathrm{const}_2\;
c_{\lambda,j;\mu,k}^{L(\Lambda)}\,t^\star,
$$
where $\mathrm{const}_1,\mathrm{const}_2\in q^{\frac1s\mathbb{Z}}$ depend
on $\Lambda$, $\lambda$, $\mu$, $q$, and $s$ is the cardinality of the
quotient group $P/Q$.
\end{proposition}

\medskip

Set $x=tt^\star=t^\star t$.

\begin{corollary}\label{quasicom_x}
If $\Lambda\in P_+$ and $\lambda,\mu$ are weights of $L(\Lambda)$ then
\begin{equation}\label{x_c}
x\,c_{\lambda,j;\mu,k}^{L(\Lambda)}=\mathrm{const}\;
c_{\lambda,j;\mu,k}^{L(\Lambda)}\,x,
\end{equation}
with $\mathrm{const}$ depending on $\Lambda$, $\lambda$, $\mu$,
$q$, and belonging to $q^{\frac1s\mathbb{Z}}$.
\end{corollary}

The multiplicative system $x^{\mathbb{Z}_+}$ is an Ore set, see \eqref{x_c}. Let
$\mathbb{C}[G]_{q,x}$ be a localization of $\mathbb{C}[G]_q$ with respect to this
multiplicative system. The fact that $\mathbb C[G]_q$ is a (left and right) Noetherian
integral domain, see \cite[p. 266]{Bou_spectr}, implies that
$\mathbb{C}[G]_q\hookrightarrow\mathbb{C}[G]_{q,x}$.

\begin{proposition}\label{*_local_gen}
\begin{itemize}
\item[1.] The $U_q\mathfrak{g}$-module algebra structure of $\mathbb{C}[G]_q$ admits a
unique extension to $\mathbb{C}[G]_{q,x}$.

\item[2.] For all $\xi\in U_q\mathfrak{g}$ and $f\in\mathbb{C}[G]_{q,x}$ there exists a unique Laurent polynomial
$p_{f,\xi}$ in one variable with coefficients from $\mathbb{C}[G]_{q,x}$such that
\begin{equation*}
\xi(f\cdot x^n)=p_{f,\xi}(q^\frac{n}s)\,x^n
\end{equation*}
for all $n\in\mathbb{Z}$.
\end{itemize}
\end{proposition}

{\bf Proof.} Consider first the second claim. The fact that $\mathbb{C}[G]_q$ is a
$U_q\mathfrak{g}$-module algebra allows one to reduce the general case $\xi\in
U_q\mathfrak{g}$ to the special case $\xi\in\{K_j^{\pm 1},E_j,F_j|\:j=1,2,\ldots,l\}$.
Within this special case, the existence of a Laurent polynomial $p_{f,\xi}$ follows from
Corollary \ref{quasicom_x} and the explicit formulas for $\triangle(K_j^{\pm 1})$,
$\triangle(E_j)$, $\triangle(F_j)$. The rest of the statements of Proposition
\ref{*_local_gen} can be proved in the same way as the corresponding claims in
Proposition \ref{*_local}. \hfill $\square$

\medskip

We thus get a well known statement.

\begin{corollary}
The multiplicative system $t^{\mathbb{Z}_+}$ is an Ore set in
$\mathbb{C}[G]_q$, so that the latter algebra is embeddable into the
localization $\mathbb{C}[G]_{q,t}$, and the structure of
$U_q\mathfrak{g}$-module algebra is uniquely extendable from
$\mathbb{C}[G]_q$ to $\mathbb{C}[G]_{q,t}$.
\end{corollary}

\begin{remark}\label{Rosenberg}
A different approach to proving that the structure of
$U_q\mathfrak{g}$-module algebra admits an extension to a localization is
presented in works by V. Lunts and A. Rosenberg \cite{LuntsRosLoc},
\cite{LuntsRosDiff2}.
\end{remark}

Now introduce the matrix element
\begin{equation*}
t'=c_{\overline{\omega}_{l_1},-\overline{\omega}_{l_0}+\alpha_{l_0}}^
{\overline{\omega}_{l_1}}
\end{equation*}
of the fundamental representation $\pi_{l_1}$ of $U_q\mathfrak{g}$ in a special basis of
the vector space $L(\overline{\omega}_{l_1})$. Note that the weights
$-\overline{\omega}_{l_0}$, $-\overline{\omega}_{l_0}+\alpha_{l_0}$ belong to the
$W$-orbit of the highest weight $\overline{\omega}_{l_1}$, hence $t$, $t'$ do not depend
on a choice of a special basis.

It follows from the definitions that for all $j\in\{1,2,\ldots,l\}$

\begin{gather}
F_j\,t=0,\qquad K^{\pm 1}_j\,t=
\begin{cases}
q_j^{\mp 1}\,t, & j=l_0,
\\ 0, & j\ne l_0,
\end{cases} \qquad
 E_j\,t=
\begin{cases}
t', & j=l_0,
\\ 0, & j\ne l_0.
\end{cases}\label{FHE_t}
\end{gather}

\begin{proposition}\label{hol_embed}
The map $i:z_\mathrm{low}\mapsto t^{-1}t'$ is uniquely extendable to an
embedding of $U_q\mathfrak{g}$-module algebras
\begin{equation}\label{hom_i}
i:\mathbb{C}[\mathfrak{p}^-]_q\hookrightarrow\mathbb{C}[G]_{q,x}.
\end{equation}
\end{proposition}

{\bf Proof.} The pair $(\mathfrak{g},\mathfrak{k})$ in question is Hermitian symmetric.
Hence the uniqueness of the homomorphism \eqref{hom_i} is due to the fact that the
$U_q\mathfrak{k}$-module $\mathbb{C}[\mathfrak{p}^-]_{q,1}$ is simple together with the
fact that the subspace $\mathbb{C}[\mathfrak{p}^-]_{q,1}$ generates the unital algebra
$\mathbb{C}[\mathfrak{p}^-]_q$. Both statements are proved in \cite{SSV4}.

We now turn to proving the existence of morphism \eqref{hom_i}. Consider
the $U_q\mathfrak{g}$-module $\widetilde{w}_0U_q\mathfrak{g}$, the dual
$U_q\mathfrak{g}$-module $(\widetilde{w}_0U_q\mathfrak{g})^*$, and its
largest weight submodule
$\mathbb{F}\subset(\widetilde{w}_0U_q\mathfrak{g})^*$ formed by
$U_q\mathfrak{b}^-$-finite elements. Here $U_q\mathfrak{b}^-$ is the
standard notation for the Hopf subalgebra generated by
$\{K_j^\pm,F_j\}_{j=1,2,\cdots,l}$. That is,
$\mathbb{F}=\bigoplus\limits_{\lambda\in P}\mathbb{F}_\lambda$, with
$$
\mathbb{F}_\lambda=\left\{\left.f\in(\widetilde{w}_0U_q\mathfrak{g})^*
\right|\:\dim(U_q\mathfrak{b}^-f)<\infty,\quad H_if=\lambda(H_i)f,\;
i=1,\ldots,l\right\}.
$$
The vector subspace $\mathbb{C}[G]_q$ is dense in $\mathbb{F}$
with respect to the weak topology. An application of
\eqref{w0_R_new} makes it easy to prove that the structure of
$U_q\mathfrak{g}$-module algebra on $\mathbb{C}[G]_q$ admits an
extension by continuity to $\mathbb{F}$. In fact,
\begin{eqnarray*}
\langle f_1f_2,\widetilde{w}_0\xi\rangle &=& \langle
R_{\mathbb{F}\mathbb{F}}\triangle(\xi)(f_1\otimes
f_2),\widetilde{w}_0\otimes\widetilde{w}_0\rangle,
\\ \langle\xi f,\widetilde{w}_0\eta\rangle &=& \langle
f,\widetilde{w}_0\eta\xi\rangle,\qquad\xi,\eta\in U_q\mathfrak{g},\qquad
f,f_1,f_2\in\mathbb{F}.
\end{eqnarray*}
Here $R_{\mathbb{F}\mathbb{F}}$ is the linear map in
$\mathbb{F}\otimes\mathbb{F}$ determined by the action of the universal
R-matrix.

It was noted in the proof of the uniqueness of the morphism \eqref{hom_i} that
$z_\mathrm{low}$ generates the $U_q\mathfrak{k}$-module algebra
$\mathbb{C}[\mathfrak{p}^-]_q$. Thus, the existence of the embedding \eqref{hom_i} is
proved if one produces embeddings of $U_q\mathfrak{g}$-module algebras
$$
i_1:\mathbb{C}[G]_{q,t}\hookrightarrow\mathbb{F},\qquad
i_2:\mathbb{C}[\mathfrak{p}^-]_q\hookrightarrow\mathbb{F},
$$
such that $i_1(t^{-1}t')=i_2(z_{\operatorname{low}})$.

To produce an embedding $i_1$, extend the natural pairing
$\mathbb{C}[G]_q\times\widetilde{w}_0U_q\mathfrak{g}\to\mathbb{C}$ to a
pairing $\mathbb{C}[G]_{q,t}\times\widetilde{w}_0U_q\mathfrak{g}\to
\mathbb{C}$. In the classical case $q=1$ such extension is possible because
$w_0$ belongs to the domain $t\ne 0$ where all the functions from
$\mathbb{C}[G]_{q,t}$ are regular. In the quantum case we are going to use
\begin{equation}\label{w0_t}
\langle t,\widetilde{w}_0\rangle\ne 0,
\end{equation}
which follows from \eqref{w_action}.

Let $f\in\mathbb{C}[G]_{q,\lambda}$ be a weight vector with weight
$\lambda\in P$ and set $c_k(f)=\langle ft^k,\widetilde{w}_0\rangle$. An
application of the equations
$$
\left\langle ft^k,\widetilde{w}_0\right\rangle=\left\langle
R_{\mathbb{C}[G]_q\mathbb{C}[G]_q}(ft^{k-1}\otimes
t),\widetilde{w}_0\otimes\widetilde{w}_0\right\rangle,
$$
\eqref{FHE_t} and an explicit formula for the universal R-matrix makes it easy to prove
that the sequence of numbers $\{c_k(f)\}$, $k\in\mathbb{Z}_+$, is a solution of a first
order difference equation whose coefficients do not depend on $f$. Specifically,
$$
c_k(f)=
q^{\mathrm{const}_1\cdot(k-1)+\mathrm{const}_2\cdot\lambda}\cdot\langle
t,\widetilde{w}_0\rangle c_{k-1}(t),\qquad k\in\mathbb{N}.
$$
Using the last formula and \eqref{w0_t}, we can introduce $\{c_k(f)\}$ for all
$k\in\mathbb{Z}$. Thus we get a `natural' extension of the linear
functional $f\mapsto\langle f,\widetilde{w}_0\rangle$ from
$\mathbb{C}[G]_q$ to $\mathbb{C}[G]_{q,t}$.

Proposition \ref{*_local_gen} implies that the subalgebra
$\mathbb{C}[G]_{q,t}\subset\mathbb{C}[G]_{q,x}$ is a
$U_q\mathfrak{g}$-module algebra. Consider the pairing
\begin{equation*}
\mathbb{C}[G]_{q,t}\times\widetilde{w}_0U_q\mathfrak{g}\to\mathbb{C},\qquad
f\times\widetilde{w}_0\xi\mapsto\langle\xi f,\widetilde{w}_0\rangle
\end{equation*}
together with the associated morphism of  $U_q\mathfrak{g}$-modules
$$
i_1:\mathbb{C}[G]_{q,t}\to(\widetilde{w}_0U_q\mathfrak{g})^*,\qquad
i_1:f\mapsto\langle f,\cdot\rangle.
$$
Let us prove that $i_1$ is injective. It follows from invertibility of
$\widetilde{w}_0\in\mathbb{C}[G]_q^*$ that the vector subspace
$\widetilde{w}_0U_q\mathfrak{g}$ is dense in $\mathbb{C}[G]_q^*$ with
respect to the weak topology. Hence
$\operatorname{Ker}i_1\cap\mathbb{C}[G]_q=0$. What remains is to note that
$i_1(f)=0$ implies $i_1(f\cdot t^j)=0$ for all $j\in\mathbb{N}$ since
$$
\langle f\cdot t^j,\widetilde{w}_0\xi\rangle=\left\langle
R_{\mathbb{C}[G]_{q,t}\mathbb{C}[G]_{q,t}}\triangle(\xi)(f\otimes
t^j),\widetilde{w}_0\otimes\widetilde{w}_0\right\rangle.
$$

Note that $i_1\mathbb{C}[G]_{q,t}\subset\mathbb{F}$ since $i_1$ is a
morphism of $U_q\mathfrak{g}$-modules and all $f\in\mathbb{C}[G]_{q,t}$ are
$U_q\mathfrak{b}^-$-finite.\footnote{It suffices to prove the last
statement for $f\in\mathbb{C}[G]_q$ and $f=t^{-1}$. On the other hand,
$\mathbb{C}[G]_q=\bigoplus\limits_{\lambda\in
P_+}\operatorname{End}L(\lambda)$, and $\dim(U_q\mathfrak{b}^-\cdot
t^{-1})=1$, as one can see from \eqref{FHE_t}}

We now prove that the embedding $i_1:\mathbb{C}[G]_{q,t}\hookrightarrow\mathbb{F}$
we have just obtained is a morphism of algebras. It suffices to show that,
given a pair of weight vectors $f_1,f_2\in\mathbb{C}[G]_q$, the equation
\begin{equation}\label{eq_f1f2}
\langle(f_1\cdot t^k)(f_2\cdot t^k),\widetilde{w}_0\rangle=\left\langle
R_{\mathbb{C}[G]_{q,t}\mathbb{C}[G]_{q,t}}(f_1\cdot t^k\otimes f_2\cdot
t^k),\widetilde{w}_0\otimes\widetilde{w}_0\right\rangle
\end{equation}
holds for all $k\in\mathbb{Z}$. As the sequences of numbers in either
sides of \eqref{eq_f1f2} are linear combinations of sequences of the form
$q^{ak^2+bk}$, $a,b\in\mathbb{Q}$, it suffices to prove \eqref{eq_f1f2} for
large $k$. For example, one may restrict considerations to the special case
$f_1,f_2\in\mathbb{C}[G]_q$. In this case \eqref{eq_f1f2} follows from
\eqref{w0_R_new}.

The embedding $i_1$ is thus produced. Now for the construction of $i_2$: let
$i_2:\mathbb{C}[\mathfrak{p}^-]_q\to(\widetilde{w}_0U_q\mathfrak{g})^*$ be
a linear map dual to
$$
\mathscr{I}:\widetilde{w}_0U_q\mathfrak{g}\to
N(\mathfrak{q}^+,0),\qquad\mathscr{I}:\widetilde{w}_0\xi\mapsto
S(\xi)v(\mathfrak{q}^+,0).
$$
It follows from the definitions that $i_2$ is an injective morphism of
$U_q\mathfrak{g}$-modules. Hence
$i_2:\mathbb{C}[\mathfrak{p}^-]_q\hookrightarrow\mathbb{F}$ since the
$U_q\mathfrak{g}$-module $\mathbb{C}[\mathfrak{p}^-]_q$ is locally
$U_q\mathfrak{b}^-$-finite dimensional.

Now properties of the universal R-matrix and the fact that $v(\mathfrak{q}^+,0)$ is the
highest vector with zero weight imply that
\begin{multline*}
\langle\xi(f_1f_2),v(\mathfrak{q}^+,0)\rangle=
\langle\triangle(\xi)(f_1\otimes f_2),v(\mathfrak{q}^+,0)\otimes
v(\mathfrak{q}^+,0)\rangle=
\\ =\langle R_{\mathbb{C}[\mathfrak{p}^-]_q\mathbb{C}[\mathfrak{p}^-]_q}
\triangle(\xi)(f_1\otimes f_2),v(\mathfrak{q}^+,0)\otimes
v(\mathfrak{q}^+,0)\rangle=
\\ =\langle
R_{\mathbb{F}\mathbb{F}}\triangle(\xi)(i_2(f_1)\otimes
i_2(f_2)),\widetilde{w}\otimes\widetilde{w}\rangle=
\langle\xi(i_2(f_1)i_2(f_2),\widetilde{w}\rangle,
\end{multline*}
with $f_1,f_2\in\mathbb{C}[\mathfrak{p}^-]_q$, $\xi\in U_q\mathfrak{g}$.
Hence $i_2$ is a morphism of algebras.

We now prove that $i_1(t^{-1}t')\in i_2\mathbb{C}[\mathfrak{p}^-]_q$. It suffices to
demonstrate that $i_1(t^{-1}t')$ is orthogonal to the subspace
$\operatorname{Ker}\mathscr{I}=\widetilde{w}_0(U_q\mathfrak{q}^+)_0$, where
\hbox{$(U_q\mathfrak{q}^+)_0= U_q\mathfrak{q}^+\cap\operatorname{Ker}\varepsilon$}, and
$\varepsilon$ is a counit of the Hopf algebra $U_q\mathfrak{q}^+$. Consider the `regular'
representation $L_{\operatorname{reg}}$ of the algebra $U_q\mathfrak{g}^{\rm op}$ in the
space $\mathbb{C}[G]_q$ defined as follows:
\begin{equation*}
L_{\operatorname{reg}}(\xi):\;f(\eta)\;\mapsto\;f(\xi\eta),
\end{equation*}
with $f\in\mathbb{C}[G]_q$, $\xi,\eta\in U_q\mathfrak{g}$. Applying the ideas of the
proof of Proposition \ref{*_local_gen}, we may extend $L_\mathrm{reg}(\xi)$ to
$\mathbb{C}[G]_{q,t}$. In this way, $\mathbb{C}[G]_{q,t}$ is furnished with a structure
of $U_q\mathfrak{g}^\mathrm{op}$-module. Given $\xi \in (U_q\mathfrak{q}^+)_0$, the
element $\eta= \widetilde{w}_0 \cdot \xi \cdot (\widetilde{w}_0)^{-1}$ of
$(\mathbb{C}[G]_q)^*$ belongs to the subalgebra $U_q\mathfrak{g,}$ and
\begin{equation}\label{last_eq_embed}
L_\mathrm{reg}(\eta)(t^{-1}t')=0.
\end{equation}
In fact, while proving this relation, $\xi$ may be assumed homogeneous of
positive degree. Equip the $U_q\mathfrak{g}^{op}$-module
$V=\{L_\mathrm{reg}(\zeta)(t^{-1}t')|\:\zeta\in U_q\mathfrak{g}^{op}\}$
with the grading corresponding to the subset $\mathbb{S}_1=\{1,2,\ldots,
l\}\setminus\{l_1\}$. It is easy to demonstrate that the element
$L_{\operatorname{reg}}(\eta) (t^{-1}t') \in V$ has incompatible
homogeneity properties, hence is zero.

Now \eqref{last_eq_embed} implies the required orthogonality.

\medskip

To finish the proof of Proposition \ref{hol_embed}, one has to
establish that the element $z$ with the property
$i_1(t^{-1}t')=i_2(z)$ is just $z_\mathrm{low}$. Recall that
$\mathbb{S}=\{1,2,\ldots,l\}\setminus\{l_0\}$ and
$\mathbb{C}[G]_{q,t}$ is a weight $U_q\mathfrak{g}$-module. Hence
we have a well defined linear map $H_\mathbb{S}$ in
$\mathbb{C}[G]_{q,t}$. The relations
\begin{eqnarray*}
H_\mathbb{S}t &=& -\overline{\omega}_{l_0}(H_\mathbb{S})t,\qquad
F_jt=0,\quad j=1,2,\ldots,l,
\\ H_\mathbb{S}t' &=&
(-\overline{\omega}_{l_0}(H_\mathbb{S})+\alpha_{l_0}(H_\mathbb{S}))t',
\qquad F_jt=0,\quad j\ne l_0,
\end{eqnarray*}
imply that
$$
H_\mathbb{S}(t^{-1}t')=2t^{-1}t',\qquad F_j(t^{-1}t')=0,\quad j\ne l_0.
$$
Thus the element $z\in\mathbb{C}[\mathfrak{p}^-]_q$ with the property
$i_1(t^{-1}t')=i_2(z)$ is the lowest weight vector of the
$U_q\mathfrak{k}$-module $\mathbb{C}[\mathfrak{p}^-]_{q,1}$. A
straightforward computation shows that
\begin{equation*}
t'=q_{l_0}^{\frac{1}{2}}E_{l_0} t,
\end{equation*}
$F_{l_0}(t^{-1}t')=q_{l_0}^{\frac{1}{2}}$. Hence
$F_{l_0}(z)=q_{l_0}^{1/2}$, and the relation
$z=z_{\operatorname{low}}$ now follows from Lemma \ref{z1}. \hfill
$\square$

\bigskip

The embedding $i:\mathbb{C}[\mathfrak{p}^-]_q\hookrightarrow
\mathbb{C}[G]_{q,x}$ just produced will be called the canonical embedding.

\begin{proposition}\label{image_hol}
The image of $\mathbb{C}[\mathfrak{p}^-]_q$ under the canonical embedding
into $\mathbb{C}[G]_{q,x}$ is the subalgebra generated by
$\{t^{-1}\,c_{\overline{\omega}_{l_1},\lambda}^{\overline{\omega}_{l_1}}\}$,
with $c_{\overline{\omega}_{l_1},\lambda}^{\overline{\omega}_{l_1}}$ being
the elements of the 'top row' of the matrix of the fundamental
representation $\pi_{l_1}$.
\end{proposition}

{\bf Proof.} Let $\mathcal{F} \subset \mathbb{C}[G]_{q,x}$ be the
subalgebra generated by $\{ t^{-1}c_{\overline{\omega}_{l_1},\lambda}
^{\overline{\omega}_{l_1}}\}$. It follows from the definitions that it is a
$U_q\mathfrak{g}$-module subalgebra.

We now establish that $i\mathbb{C}[\mathfrak{p}^-]_q=\mathcal{F}$. To begin
with, observe that
$i\mathbb{C}[\mathfrak{p}^-]_q\subset\mathcal{F}$ since
$i(z_\mathrm{low})\in\mathcal{F}$, hence
$i\mathbb{C}[\mathfrak{p}^-]_{q,1} \subset\mathcal{F}$. The
opposite inclusion
$i\mathbb{C}[\mathfrak{p}^-]_q\supset\mathcal{F}$ is easily
accessible via an application of the fact that the linear span
$\mathcal{L}$ of the set
$\{c_{\overline{\omega}_{l_1},\lambda}^{\overline{\omega}_{l_1}}\}$
is just $U_q\mathfrak{b}^+\cdot t$, since $t$ is the lowest weight
vector of the $U_q\mathfrak{g}$-module $\mathcal{L} \approx
L(\overline{\omega}_{l_0})$. \hfill $\square$

\begin{remark}\label{mult_free}
An application of the canonical embedding, together with the well-known Hua-Schmid
theorem (see., e.g., \cite[p. 73]{John}, \cite[p. 443]{Takeuchi}), allows one easily to
establish that
$$
\dim\operatorname{Hom}_{U_q\mathfrak{k}}
(L(\mathfrak{k},\lambda),L(j\overline{\omega}_{l_0}))\le 1
$$
for all $\lambda\in P_+^\mathbb{S}$, $j\in\mathbb{N}$.
\end{remark}

\bigskip

Let $\mathbb{C}[X^-_{\mathbb{S}}]_q\subset\mathbb{C}[G]_q$ be the smallest
$U_q\mathfrak{g}$-module subalgebra which contains $t$, and
$\mathbb{C}[X^+_{\mathbb{S}}]_q\subset\mathbb{C}[G]_q$ the smallest
$U_q\mathfrak{g}$-module subalgebra which contains $t^\star$. This can also be
formulated as follows: the subalgebra $\mathbb{C}[X^-_{\mathbb{S}}]_q$ is
generated by the matrix elements
$\left\{c_{\overline{\omega}_{l_1};\lambda,j}^{\overline{\omega}_{l_1}}
\right\}$, and the subalgebra $\mathbb{C}[X^+_{\mathbb{S}}]_q$ is generated
by
$\left\{c_{w_0\overline{\omega}_{l_0};\lambda,j}^{\overline{\omega}_{l_0}}
\right\}$, see \eqref{aster_standard}.

We equip the algebra $\mathbb{C}[X^-_{\mathbb{S}}]_q$ with a grading by
setting
\hbox{$\deg\left(c_{\overline{\omega}_{l_1};\lambda,j}^{\overline{\omega}_{l_1}}
\right)=1$.} This grading is well-defined since it can be imposed
in a different way by using the linear map
$L_\mathrm{reg}(H_{\mathbb{S}_1})$. Specifically, $\deg f=j$ if and only if
$L_\mathrm{reg}(H_{\mathbb{S}_1})f=jf$.

\medskip

The above results imply the following statements, which are well-known in
quantum group theory.

Firstly, the multiplicative system $t^{\mathbb{Z}_+}$ is an Ore set in
$\mathbb{C}[X^-_\mathbb{S}]_q$ , and this algebra is naturally embeddable into its
localization: $\mathbb{C}[X^-_\mathbb{S}]_q
\hookrightarrow\mathbb{C}[X^-_\mathbb{S}]_{q,t}$.

Secondly, $\mathbb{C}[X^-_{\mathbb{S}}]_q$ is a
$U_q\mathfrak{g}$-module subalgebra.

Thirdly, the structure of $U_q\mathfrak{g}$-module algebra on
$\mathbb{C}[X^-_{\mathbb{S}}]_q$ admits an extension to
$\mathbb{C}[X^-_{\mathbb{S}}]_{q,t}$, and this extension is unique.

\medskip

Now Proposition \ref{hol_embed} can be used to obtain

\begin{corollary}\label{embed_to_flags}
The map $i:z_\mathrm{low}\mapsto t^{-1}t'$ is uniquely extendable to an
embedding of $U_q\mathfrak{g}$-module algebras
$i:\mathbb{C}[\mathfrak{p}^-]_q\hookrightarrow
\mathbb{C}[X^-_{\mathbb{S}}]_{q,t}$.
\end{corollary}

\begin{remark}\label{embed-star}
In a similar way one can equip the localization
$\mathbb{C}[X_\mathbb{S}^+]_{q,t^\star}$ of
$\mathbb{C}[X_\mathbb{S}^+]_q$ with respect to the multiplicative
system $(t^\star)^{\mathbb{Z}_+}$ with the structure of a
$U_q\mathfrak{g}$-module algebra and obtain an embedding
$\mathbb{C}[\mathfrak{p}^+]_q\hookrightarrow
\mathbb{C}[X_\mathbb{S}^+]_{q,t^\star}$.
\end{remark}

Consider the classical case $q=1$. We are about to use an isomorphism of
$U\mathfrak{g}$-modules $L(\overline{\omega}_{l_1})^*\approx L(\overline{\omega}_{l_0})$.
Consider the cone $X^-_\mathbb{S}\subset L(\overline{\omega}_{l_0})$ generated by the
$G$-orbit of the highest weight vector $v(\overline{\omega}_{l_0})\in
L(\overline{\omega}_{l_0})$. It is a homogeneous space of the group $G$, and the
stability group of the point $\mathbb{C}v(\overline{\omega}_{l_0})$ is a standard
parabolic subgroup \cite[p. 291]{Hum}.

We have produced a quantum analog of the embedding of $\mathfrak{p}^-$ into the above
homogeneous space, which is called a generalized flag variety. In the classical case
$q=1$ such an embedding is well-known and has been extensively investigated \cite[p. 34,
35]{Springer}. In the special case $\mathfrak{g}=\mathfrak{sl}_N$ a result close to ours
has been obtained by Fioresi \cite{Fioresi_cell}.


\medskip

\subsection{A $*$-algebra \boldmath
$(\mathbb{C}[\mathbb{X}_\mathbb{S}]_{q},*)$.}\label{star_emb}

Consider the smallest unital $U_q\mathfrak{g}$-module subalgebra
$\mathbb{C}[\mathbb{X}_\mathbb{S}]_q\subset\mathbb{C}[G]_q$ which
contains $t$, $t^\star$.

It is known \cite[p. 480]{Stokman_list} that
$\mathbb{C}[\mathbb{X}_\mathbb{S}]_q$ consists of functions
$f\in\mathbb{C}[G]_q$ such that
$$
L_\mathrm{reg}(E_j)f=L_\mathrm{reg}(F_j)f=L_\mathrm{reg}(K_j^{\pm
1}-1)f=0, \qquad j\neq l_1.
$$

\begin{example} If $\mathfrak{g}=\mathfrak{sl}_2$, then
$$\mathbb{C}[\mathbb{X}_\mathbb{S}]_q=\mathbb{C}[SL_2]_q.$$
\end{example}

\begin{proposition} \label{star_X_new}
There exists a unique antilinear involution $*$ in
$\mathbb{C}[\mathbb{X}_\mathbb{S}]_q$ such that
$(\mathbb{C}[\mathbb{X}_\mathbb{S}]_q,*)$ is a
$(U_q\mathfrak{g},*)$-module algebra and $t^*=t^\star$.
\end{proposition}

{\bf Proof.} The uniqueness of the required involution is obvious.
We prove its existence by equipping $\mathbb{C}[G]_q$ with a
$(U_q\mathfrak{g},*)$-module algebra structure such that $t^*=t^\star$.

Firstly, we apply the following natural embeddings
$U_q\mathfrak{g}\hookrightarrow\mathbb{C}[G]_q^*$,
$\widetilde{w}_0U_q\mathfrak{g}\hookrightarrow\mathbb{C}[G]_q^*$. The image
of $\widetilde{w}_0U_q\mathfrak{g}$ under the embedding is dense in
\hbox{$\mathbb{C}[G]_q^*\cong\mathop{\times}\limits_{\lambda \in P_+}
\operatorname{End}L(\lambda)$} in the product topology since
$\widetilde{w}_0$ is invertible. Therefore the canonical pairing
$\langle\cdot,\cdot\rangle$
$$\mathbb{C}[G]_q\times\widetilde{w}_0U_q\mathfrak{g}\to\mathbb{C}$$
is non-degenerate. Now we need some auxiliary observations.

\begin{lemma}\label{def_star}
There exists a unique antilinear map $*$ in $\mathbb{C}[G]_q$ such that
\begin{equation}\label{dual_star} \langle
f^*,\widetilde{w}_0\xi\rangle=\overline{\langle
f,\widetilde{w}_0(S(\xi))^*\rangle}
\end{equation}
for any $f\in\mathbb{C}[G]_q$, $\xi\in U_q\mathfrak{g}$.
\end{lemma}

{\bf Proof.} The uniqueness is obvious. We establish the existence by using the
involution $\sharp$ in $\mathbb{C}[G_{\mathbb R}]_q$ (see \ref{canon_example}).

The representation $L_{\rm reg}$ admits an extension by continuity to a representation
$\overline{L}_\mathrm{reg}$ of $\mathbb{C}[G]_q^{*\mathrm{op}}$. Comparing
\eqref{dual_star} with \eqref{g_0}, one has that for all $f\in\mathbb{C}[G]_q$, $\xi\in
U_q\mathfrak{g}$
$$
\left\langle\overline{L}_\mathrm{reg}(\widetilde{w}_0)f^*,\xi\right\rangle=
\left\langle\overline{L}_\mathrm{reg}(\widetilde{w}_0)f,
(S(\xi))^*\right\rangle=
\left\langle\left(\overline{L}_\mathrm{reg}(\widetilde{w}_0)f\right)^\sharp,
\xi\right\rangle.
$$
Hence,
\begin{equation}\label{intertw_stars}
\overline{L}_\mathrm{reg}(\widetilde{w}_0)\cdot
*=\sharp\cdot\overline{L}_\mathrm{reg}(\widetilde{w}_0).
\end{equation}
Finally,
\begin{equation*}
f^*=L_\mathrm{reg}(\widetilde{w}_0^{-1})
\left(L_\mathrm{reg}(\widetilde{w}_0)f\right)^\sharp.
\end{equation*}
The last equation yields an anti-linear operator with the
required properties. \hfill $\square$

\bigskip

We now prove that $*$ is an involution of $\mathbb{C}[G]_q$.
\begin{lemma}\label{prop_6.1}
1. $f^{**}=f$ for any $f\in\mathbb{C}[G]_q$.

2. $(f_1f_2)^*=f_2^*f_1^*$ for any $f_1,f_2\in\mathbb{C}[G]_q$.
\end{lemma}

{\bf Proof.} The square of the antilinear operator
$*S:U_q\mathfrak{g}\to U_q\mathfrak{g}$ is the identity map. This
proves the first statement:
$$
\langle f^{**},\widetilde{w}_0\xi\rangle=\langle
f,\widetilde{w}_0(S((S(\xi))^*))^*\rangle=\langle f,\xi\rangle
$$
for any $f\in\mathbb{C}[G]_q$, $\xi\in U_q\mathfrak{g}$.

Consider the second statement. Strictly speaking, the universal R-matrix
should be substituted by the corresponding map
$R_{\mathbb{C}[G]_q\mathbb{C}[G]_q}$ in the coming proof. For the sake of
clarity, we do not do this. We use some properties of the operator $*S$.
These properties come from the following equalities
\begin{gather*}
\triangle(S(\xi))=(S\otimes
S)\triangle^\mathrm{cop}(\xi),\qquad\triangle(\xi^*)= (\triangle(\xi))^{*\otimes
*},
\\ S\otimes S(R)=R,\qquad R^{*\otimes *}=R_{21},
\end{gather*}
where $R_{21}$ is obtained from $R$ via a flip of tensor factors.
The last equality follows from the well-known properties of the
universal R-matrix, see, for example Proposition 2.3.1 in
\cite{KorSoib}.

For any $f_1,f_2\in\mathbb{C}[G]_q$, $\xi\in U_q\mathfrak{g}$
$$
\langle(f_1f_2)^*,\widetilde{w}_0\xi\rangle=\overline{\langle
f_1f_2,\widetilde{w}_0(S(\xi))^*\rangle}=\overline{\langle
f_1\otimes
f_2,(\widetilde{w}_0\otimes\widetilde{w}_0)R\triangle(S(\xi))^*\rangle};
$$
\begin{multline*}
\langle f_2^*f_1^*,\widetilde{w}_0\xi\rangle=\langle f_2^*\otimes
f_1^*,(\widetilde{w}_0\otimes\widetilde{w}_0)R\triangle(\xi)\rangle=
\\ =\overline{\langle f_2\otimes
f_1,(\widetilde{w}_0\otimes\widetilde{w}_0)(*S\otimes
*S)(R\triangle(\xi))\rangle}=
\\ =\overline{\langle f_2\otimes
f_1,(\widetilde{w}_0\otimes\widetilde{w}_0))
R_{21}\triangle^\mathrm{cop}((S(\xi))^*)\rangle}=
\\ =\overline{\langle f_1\otimes
f_2,(\widetilde{w}_0\otimes\widetilde{w}_0))R\triangle((S(\xi))^*)\rangle}.
\end{multline*}
Therefore $\langle(f_1f_2)^*,\widetilde{w}_0\xi\rangle=\langle
f_2^*f_1^*,\widetilde{w}_0\xi\rangle$ for any
$f_1,f_2\in\mathbb{C}[G]_q$, $\xi\in U_q\mathfrak{g}$. \hfill
$\square$

\bigskip

Consider the $*$-algebra $\mathbb{C}[w_0G_{\mathbb R
}]_q\overset{\mathrm{def}}{=}(\mathbb{C}[G]_q,*)$. It is a q-analog of the algebra
$\mathbb{C}[w_0G_{\mathbb R}]$ of regular functions on the real affine algebraic variety
$w_0G_{\mathbb R}$, where $G_{\mathbb R}$ is a noncompact real form of $G$.

\begin{example}
If $\mathfrak{g}=\mathfrak{sl}_2$, then $(\alpha,\alpha)=2$,
$(\overline{\omega},\overline{\omega})=\frac12$, $(H,H)=2$. Using
\eqref{s_bar_eq}, one can explicitly find the action of $\widetilde{w}_0$ in the
fundamental representation $\overline{\pi}$ of the algebra
$\mathbb{C}[SL_2]_q^*$:
$$
\overline{\pi}(\widetilde{w}_0)=q^{-1/4}\begin{pmatrix}0 & -1\\
q^{-1} & 0\end{pmatrix}.
$$
Hence,
\begin{equation*}
\begin{pmatrix}\overline{L}_\mathrm{reg}(\widetilde{w}_0)t_{11} &\;
\overline{L}_\mathrm{reg}(\widetilde{w}_0)t_{12}\\
\overline{L}_\mathrm{reg}(\widetilde{w}_0)t_{21}&\;
\overline{L}_\mathrm{reg}(\widetilde{w}_0)t_{22}\end{pmatrix}=
q^{-1/4}\begin{pmatrix}0 &\; -1\\ q^{-1} &\;
0\end{pmatrix}\begin{pmatrix}t_{11} &\; t_{12}\\ t_{21} &\;
t_{22}\end{pmatrix}.
\end{equation*}

One can easily check that the involutions \eqref{inv},
\eqref{inv_su11} satisfy \eqref{intertw_stars}. For example,
 $$
 \overline{L}_\mathrm{reg}(\widetilde{w}_0)t_{11}^*=
\overline{L}_\mathrm{reg}(\widetilde{w}_0)(-t_{22})=-q^{-5/4}t_{12}
=(-t_{21})^\sharp=
\left(\overline{L}_\mathrm{reg}(\widetilde{w}_0)t_{11}\right)^\sharp.
$$
Hence $\overline{L}_\mathrm{reg}(\widetilde{w}_0)f^*=
\left(\overline{L}_\mathrm{reg}(\widetilde{w}_0)f\right)^\sharp$
for $f=t_{11}$.

Therefore one can see that the new definition of the involution $*$
generalizes the previous definition \eqref{inv}, which is applicable only
in the special case $\mathfrak{g}=\mathfrak{sl}_2$.
\end{example}

\begin{lemma}\label{mod_alg}
$\mathbb{C}[w_0G_{\mathbb R}]_q$ is a $(U_q\mathfrak{g},*)$-module algebra.
\end{lemma}

{\bf Proof.} For any $f\in\mathbb{C}[G]_q$, $\xi,\eta\in U_q\mathfrak{g}$ we
have that
$$
\langle(\xi f)^*,\widetilde{w}_0\eta\rangle=\overline{\langle\xi
f,\widetilde{w}_0(S(\eta))^*\rangle}=\overline{\langle
f,\widetilde{w}_0(S(\eta))^*\xi\rangle};
$$
\begin{multline*}
\langle (S(\xi))^*f^*,\widetilde{w}_0\eta\rangle=\langle
f^*,\widetilde{w}_0\eta(S(\xi))^*\rangle=\overline{\langle
f,\widetilde{w}_0(S(\eta(S(\xi))^*))^*\rangle}=
\\ =\overline{\langle
f,\widetilde{w}_0(S(\eta))^*(S((S(\xi))^*))^*\rangle}=\overline{\langle
f,\widetilde{w}_0(S(\eta))^*\xi\rangle}.
\end{multline*}
Hence for any $f\in\mathbb{C}[G]_q$, $\xi\in U_q\mathfrak{g}$ one
has
$$(\xi f)^*=(S(\xi))^*f^*.\eqno\square$$

\bigskip
We now prove that $t^*=t^\star$.

\begin{lemma}\label{up_to_const}
\begin{equation*}
t^*=\mathrm{const}\cdot
c_{-\overline{\omega}_{l_1},\overline{\omega}_{l_0}}
^{\overline{\omega}_{l_0}}.
\end{equation*}
\end{lemma}

{\bf Proof.} Consider the fundamental $U_q\mathfrak{g}$-module
$L(\overline{\omega}_{l_0})$ and the corresponding representation
$\pi_{l_0}$ together with its dual $\pi_{l_0}^*$. The elements
$t^*$, $c_{-\overline{\omega}_{l_1},\overline{\omega}_{l_0}}
^{\overline{\omega}_{l_0}}$ belong to the same weight subspace of
the simple \hbox{$U_q\mathfrak{g}^\mathrm{op}\otimes
U_q\mathfrak{g}$}-module
$\left\{f\in(U_q\mathfrak{g})^*\left|\:\operatorname{Ker}f\supset
\operatorname{Ker}\pi_{l_0}^*\right.\right\}\subset\mathbb{C}[G]_q$,
which is one dimensional. \hfill $\square$

\medskip

\begin{lemma}\label{star_aster}
$t^*=t^\star$.
\end{lemma}

{\bf Proof.} It is sufficient to prove that $t^*=
(-q)^{(\overline{\omega}_{l_1}+\overline{\omega}_{l_0},\rho)}
c_{-\overline{\omega}_{l_1},\overline{\omega}_{l_0}}
^{\overline{\omega}_{l_0}}$ since $t^\star=
(-q)^{(\overline{\omega}_{l_1}+\overline{\omega}_{l_0},\rho)}
c_{-\overline{\omega}_{l_1},\overline{\omega}_{l_0}}
^{\overline{\omega}_{l_0}}$ (see \eqref{aster_standard}).
Therefore one should establish the following
\begin{equation*} \left\langle
(c_{\overline{\omega}_{l_1},-\overline{\omega}_{l_0}}
^{\overline{\omega}_{l_1}})^*,\widetilde{w}_0\right\rangle
=\overline{\left\langle
c_{\overline{\omega}_{l_1},-\overline{\omega}_{l_0}}
^{\overline{\omega}_{l_1}},\widetilde{w}_0\right\rangle}=
(-q)^{(\overline{\omega}_{l_1}-w_0\overline{\omega}_{l_1},\rho)}\left\langle
c_{-\overline{\omega}_{l_1},\overline{\omega}_{l_0}}
^{\overline{\omega}_{l_0}},\widetilde{w}_0\right\rangle,
\end{equation*}
or, equivalently (see \eqref{w0_tilda}),
\begin{equation}\label{star_aster_2}
\overline{\left\langle
c_{\overline{\omega}_{l_1},-\overline{\omega}_{l_0}}
^{\overline{\omega}_{l_1}},\overline{w}_0^{-1}\right\rangle}=
(-q)^{(\overline{\omega}_{l_1}-
u_1\,\overline{\omega}_{l_1},\rho)} \left\langle
c_{-\overline{\omega}_{l_1},\overline{\omega}_{l_0}}
^{\overline{\omega}_{l_0}},\overline{w}_0^{-1}\right\rangle,
\end{equation}
Here $u_1$ is the shortest element in the coset
$w_0W_{\mathbb{S}_1}$ of the subgroup $W_{\mathbb{S}_1}\subset W$
generated by $\{s_j\;|\; j\neq l_1\}$.

If $\overline{\omega}_{l_1}- u_1\overline{\omega}_{l_1}=
\sum\limits_kn_k\alpha_k$, then
$(\overline{\omega}_{l_1}-u_1\overline{\omega}_{l_1},\rho)=
\sum\limits_kn_kd_k$ since
$(\alpha_k,\overline{\omega}_j)=d_j\delta_{jk}$. Hence if one
considers a reduced expression $u_1=s_{i_1}s_{i_2}\ldots s_{i_M}$
of $u_1$, one gets
$$
(-q)^{(\overline{\omega}_{l_1}-u_1 \overline{\omega}_{l_1},\rho)}=
\prod_{j=1}^M(-q_{i_j}).
$$
It follows from the last equality and the definition of
$\overline{w}_0^{\,-1}$ that it suffices to prove
\eqref{star_aster_2} in the simplest special case
$\mathfrak{g}=\mathfrak{sl}_2$. Here we can just use \eqref{s_bar}.
\hfill $\square$

\medskip

Proposition \ref{star_X_new} follows from the lemmas above. \hfill $\square$

\bigskip

In the sequel we use the notation $\mathbb{S}_1=\{1,2, \cdots, l\}
\setminus \{l_1\}$ as in the proof of Lemma \ref{star_aster}.

\begin{remark}\label{explicit_*}
One has the following equation for the matrix elements in a special basis of
$L(\overline{\omega}_{l_1})$
$$
\left(c_{\overline{\omega}_{l_1},\lambda}^{\overline{\omega}_{l_1}}\right)^*=
(-1)^{\lambda(H_{\mathbb{S}_1})+\overline{\omega}_{l_0}(H_{\mathbb{S}_1})}
\left(c_{\overline{\omega}_{l_1},\lambda}^{\overline{\omega}_{l_1}}\right)
^\star,\qquad\lambda\in W\overline{\omega}_{l_0}.
$$
In particular, for $l=m+n-1$, $l_0=n$,
$\mathfrak{g}=\mathfrak{sl}_{m+n}$, the involution $*$ can be
written explicitly \cite[p. 859]{SSV4}:
\begin{equation*}
t_{ij}^*=
\mathrm{sign}\left(\left(i-m-\frac12\right)\left(j-n-\frac12\right)\right)
t_{ij}^\star.
\end{equation*}
\end{remark}


\medskip

\subsection{The canonical embedding of \boldmath
$\operatorname{Pol}(\mathfrak{p}^-)_q$.}\label{canon_embed}

Consider the subalgebra
$\mathbb{C}[\mathbb{X}_\mathbb{S}]_{q,x}\subset\mathbb{C}[G]_{q,x}$
generated by $\mathbb{C}[\mathbb{X}_\mathbb{S}]_q$ and $x^{-1}$. It is a
localization of the algebra $\mathbb{C}[\mathbb{X}_\mathbb{S}]_q$ with
respect to the multiplicative system $x^{\mathbb{Z}_+}$.

The involution $*$ is uniquely extendable from $\mathbb{C}[\mathbb{X}_\mathbb{S}]_q$ to
$\mathbb{C}[\mathbb{X}_\mathbb{S}]_{q,x}$, and we have $(x^{-1})^*=x^{-1}$. It is easy to
prove that the structure of $U_q\mathfrak{g}$-module algebra and the involution are
compatible, i.e. we get the $(U_q\mathfrak{g},*)$-module algebra
$(\mathbb{C}[\mathbb{X}_\mathbb{S}]_{q,x},*)$, see Proposition \ref{*_local_gen}.

Similarly, one can introduce the $(U_q\mathfrak{g},*)$-module
algebra
$$\mathbb{C}[w_0G_{\mathbb R}]_{q,x}\overset{\mathrm{def}}{=}(\mathbb{C}[G]_{q,x},*).$$
\begin{theorem}\label{c_e}\begin{itemize}
\item[1.] There exists a unique homomorphism of $(U_q\mathfrak{g},*)$-module
algebras $\mathcal{I}:\operatorname{Pol}(\mathfrak{p}^-)_q \to
\mathbb{C}[\mathbb{X}_\mathbb{S}]_{q,x}$ such that $\mathcal{I}:
\,z_{\mathrm{low}}\mapsto t^{-1}t'$.
 \item[2.] $\mathcal{I}$ is injective.
\end{itemize}
\end{theorem}

{\bf Proof.} Let $\mathbb{C}[\mathbb{X}_\mathbb{S}]_{q,t}$ be the subalgebra in
$\mathbb{C}[\mathbb{X}_\mathbb{S}]_{q,x}$ generated by
$\mathbb{C}[\mathbb{X}_\mathbb{S}]_q$ and $t^{-1}$.

It follows from Proposition \ref{hol_embed} that the map
$i:z_\mathrm{low}\mapsto t^{-1}t'$ is uniquely extendable to an
embedding of $U_q\mathfrak{g}$-module algebras $
i:\mathbb{C}[\mathfrak{p}^-]_q\hookrightarrow\mathbb{C}[\mathbb{X}_\mathbb{S}]_{q,t}
$.

Furthermore, $\operatorname{Pol}(\mathfrak{p}^-)_q=$
\hbox{$(\mathbb{C}[\mathfrak{p}^-\oplus\mathfrak{p}^+],*)$}, and the
multiplication
$$m:\mathbb{C}[\mathfrak{p}^-]_q\otimes\mathbb{C}[\mathfrak{p}^+]_q\to
\mathbb{C}[\mathfrak{p}^-\oplus\mathfrak{p}^+]_q, \qquad
 m:f_-\otimes f_+\mapsto f_-f_+,$$
is a bijection.

Hence the homomorphism $\mathcal{I}$ is unique, and, if it exists,
one has
\begin{equation}\label{def_I}
\mathcal{I}(f_-f_+)=i(f_-)(i(f_+))^*,\qquad
f_\pm\in\mathbb{C}[\mathfrak{p}^\pm]_q.
\end{equation}

\medskip

Now we prove the existence. Consider the linear map
$\mathcal{I}:\operatorname{Pol}(\mathfrak{p}^-)_q\to
\mathbb{C}[\mathbb{X}_\mathbb{S}]_{q,x}$ defined by \eqref{def_I}. Here, $\mathcal{I}$ is a morphism of
$U_q\mathfrak{g}$-modules since
$\operatorname{Pol}(\mathfrak{p}^-)_q$ and
$\mathbb{C}[\mathbb{X}_\mathbb{S}]_{q,x}$ are
$(U_q\mathfrak{g},*)$-module algebras. It remains to prove that
$\mathcal{I}$ is an algebra homomorphism.

Consider the embedding
$\mathbb{C}[\mathbb{X}_\mathbb{S}]_{q,x}\subset\mathbb{C}[G]_{q,x}$. It is
sufficient to obtain the equality
\begin{equation}\label{eq_m}
(\mathcal{I}f_+)(\mathcal{I}f_-)=
m\check{R}_{\mathcal{I}\mathbb{C}[\mathfrak{p}^+]_q,
\mathcal{I}\mathbb{C}[\mathfrak{p}^-]_q}
(\mathcal{I}f_+\otimes\mathcal{I}f_-),
\end{equation}
where $f_\pm\in\mathcal{I}\mathbb{C}[\mathfrak{p}^\pm]_{q,\pm 1}$, and $m$
is the multiplication in $\mathbb{C}[G]_{q,x}$. Instead of \eqref{eq_m} we
prove a more general statement:
\begin{equation}\label{eq_gen}
q^{ck^2}\mathcal{I}(f_+)t^{*k}\cdot t^k\cdot\mathcal{I}(f_-)=
m\check{R}_{\mathcal{I}\mathbb{C}[G]_{q,t^*},\mathcal{I}\mathbb{C}[G]_{q,t}}
(\mathcal{I}(f_+)t^{*k}\otimes t^k\mathcal{I}(f_-)),
\end{equation}
where $k\in\mathbb{Z}$,\
$\mathbb{C}[G]_{q,t^*}=\{f^*|\:f\in\mathbb{C}[G]_{q,t}\}$, and $c$
is a rational number that depends on $f_-,f_+$. The last equality
can be proved under the additional assumption that
\begin{equation}\label{special_case}
\mathcal{I}(f_+)t^{*k}\in\mathbb{C}[X_\mathbb{S}^+]_q,\qquad
\mathcal{I}(f_-)t^k\in\mathbb{C}[X^-_\mathbb{S}]_q,
\end{equation}
where $\mathbb{C}[X^-_\mathbb{S}]_q$,
$\mathbb{C}[X_\mathbb{S}^+]_q\subset
\mathbb{C}[\mathbb{X}_\mathbb{S}]_q $ are $U_q\mathfrak{g}$-module
subalgebras generated by $t$, $t^*$, respectively (see also
section \ref{w0_G_x}).

Note that
$\mathbb{C}[X^+_\mathbb{S}]_q=(\mathbb{C}[X_\mathbb{S}^-]_q)^*$.

We equip $\mathbb{C}[\mathbb{X}_\mathbb{S}]_q$ with the grading via
$$
\deg\left(c_{\overline{\omega}_{l_1};\lambda,j}^{\overline{\omega}_{l_1}}
\right)=1,\qquad
\deg\left(c_{-\overline{\omega}_{l_1};\mu,k}^{\overline{\omega}_{l_0}}\right)
=-1.
$$
As in section \ref{w0_G_x}, one can prove that the grading is
well-defined.

The commutation relations in $\mathbb C[G]_q$ (see also the explicit formula for the
universal enveloping algebra and the definition of the grading in
$\mathbb{C}[\mathbb{X}_\mathbb{S}]_q$), and the equalities
$$
L_\mathrm{reg}(F_j)\psi_+=L_\mathrm{reg}(E_j)\psi_-=0,\qquad
j=1,2,\ldots,l,
$$
$$ L_\mathrm{reg}(H_j)\psi_\pm=0,\quad j\ne l_1, \qquad
(H_{\mathbb{S}},H_{\mathbb{S}})=(H_{\mathbb{S}_1},H_{\mathbb{S}_1}),
$$ imply that
\begin{equation}\label{c_r_new}
q^{-\deg(\psi_+)\deg(\psi_-)/(H_\mathbb{S},H_\mathbb{S})}\psi_+\psi_-=
m\check{R}_{\mathbb{C}[G]_q,\mathbb{C}[G]_q}
(\psi_+\otimes\psi_-),
\end{equation}
where $\psi_+\in \mathbb{C}[X_\mathbb{S}^+]_q]_q$, $\psi_-\in
\mathbb{C}[X^-_\mathbb{S}]_q$.

Using \eqref{special_case}, one can prove \eqref{eq_gen} from \eqref{c_r_new}.

The existence and uniqueness of the homomorphism $\mathcal{I}$ are now
proved. Let us prove its injectivity. We intend to prove this without
assuming $q$ to be transcendental. Our proof requires some auxiliary
results.

\bigskip

The proof consists of a construction of a representation
$\mathscr{T}$ of $\mathbb{C}[G]_{q,x}$ such that the corresponding
representation $\mathscr{T}\circ\mathcal{I}$ of
$\mathbb{C}[\mathfrak{p}^-\oplus\mathfrak{p}^+]_q$ is faithful.
This will prove the injectivity of $\mathcal{I}$.

In the construction of $\mathscr{T}$ we follow the standard approach for producing
$*$-representations of $\mathbb{C}[K]_q=(\mathbb{C}[G]_q,\star)$ which is described in
\cite{Jant}. Indeed, let $w_0^\mathbb{S}$ be the shortest element in the coset
$W_\mathbb{S}\,w_0\subset W$ of the subgroup $W_\mathbb{S}$ generated by the simple
reflections $s_i$, $i\ne l_0$, see \cite[p. 19]{Hum2}. Choose a reduced expression
$$w_0^\mathbb{S}=s_{i_1}s_{i_2}\ldots s_{i_N}$$
and introduce an irreducible representation of the Hopf $*$-algebra
$\mathbb{C}[K]_q$ in $l^2(\mathbb{Z}_+)^{\otimes N}$
$$\Pi=\Pi_{i_1}\otimes\Pi_{i_2}\otimes\cdots\otimes\Pi_{i_N},$$
where $\Pi_{i_j}$ are the $*$-representations of $(\mathbb{C}[G]_q,\star)$ in
$l^2(\mathbb{Z}_+)$ defined as follows:
$$\Pi_i=\Pi\circ\psi_i.$$
Here $\Pi$ is the $*$-representation of $\mathbb C[SU_2]_q$ defined in Remark \ref{SU_2}
and
$$\psi_i:(\mathbb{C}[G]_q,\star)\to\mathbb{C}[SU_2]_{q_i},\qquad
i=1,2,\ldots,l,$$ are the Hopf $*$-algebra homomorphisms conjugate to the
Hopf $*$-algebras embeddings
$$U_{q_i}\mathfrak{sl}_2 \hookrightarrow U_q\mathfrak{g},\qquad
K^{\pm1} \mapsto K_i^{\pm1},\quad E \mapsto E_i,\quad F \mapsto F_i.$$

It follows from the definitions of $\Pi_{i_j}$ and $\Pi$ that the
operators $\Pi(f)$, $f\in\mathbb{C}[G]_q$, admit a restriction
to the linear span $\mathscr{L}$ of elements of the standard
basis $\{e_{k_1}\otimes e_{k_2}\otimes\cdots\otimes e_{k_N}\}$ of
the Hilbert space $l^2(\mathbb{Z}_+)^{\otimes N}$. It is known
\cite[p. 314, 315]{DijStok1} that these elements of the basis are
eigenvectors of the bounded self-adjoint compact linear operator
$\Pi(x)$, and the corresponding eigenvalues are positive. Besides,
the eigenspace of the highest eigenvalue is one dimensional, and
$\mathbf{e}=e_0\otimes e_0\otimes\cdots\otimes e_0$ belongs to it.
Hence the linear operator $\Pi(x)|_\mathscr{L}$  is invertible in
the pre-Hilbert space $\mathscr{L}$. Therefore the representation
$\Pi|_\mathscr{L}$ of $\mathbb{C}[G]_q$ is canonically extendable
to the representation $\mathscr{T}$ of its localization
$\mathbb{C}[G]_{q,x}$
$$
\mathscr{T}(f)v=\Pi(f)v,\qquad f\in\mathbb{C}[G]_q,\quad
v\in\mathscr{L}.
$$
It remains to prove the faithfulness of the representation
$\mathscr{T}_\mathcal{I}=\mathscr{T}\cdot\mathcal{I}$ of
$\mathbb{C}[\mathfrak{p}^-\oplus\mathfrak{p}^+]_q$. First of all, we establish the
existence of an involution $\star$ in
$\mathbb{C}[\mathfrak{p}^-\oplus\mathfrak{p}^+]_q$ such that
\begin{equation}\label{new_*}
\mathcal{I}(f^\star)=\mathcal{I}(f)^\star,\qquad
f\in\mathbb{C}[\mathfrak{p}^-\oplus\mathfrak{p}^+]_q.
\end{equation}
Indeed, put
\begin{equation*}
f^\star\stackrel{\operatorname{def}}{=}(-1)^{\deg f}f^*,
\end{equation*}
where the grading is defined as follows:
$$ \mathbb{C}[\mathfrak{p}^-\oplus\mathfrak{p}^+]_{q,j}=
\left\{\left.f\in\mathbb{C}[\mathfrak{p}^-\oplus\mathfrak{p}^+]_q\right|\:
H_\mathbb{S}f=2jf\right\},\qquad j\in\mathbb{Z}.
$$

We introduce the notation $\operatorname{Pol}(U)_q=
(\mathbb{C}[\mathfrak{p}^-\oplus\mathfrak{p}^+]_q,\star)$. Note
that $\operatorname{Pol}(U)_q$ is a
$(U_q\mathfrak{g},\star)$-module algebra, since
$\operatorname{Pol}(\mathfrak{p}^-)_q$ is a
$(U_q\mathfrak{g},*)$-module algebra.

Now \eqref{new_*} follows from the fact that $\mathcal{I}(f^*)=\mathcal{I}(f)^*$ for any
$f\in\operatorname{Pol}(\mathfrak{p}^-)_q$ and the definitions of the involutions. Hence
$\mathscr{T}_\mathcal{I}$ is a \hbox{$*$-representation} of $\operatorname{Pol}(U)_q$ in
the pre-Hilbert space $\mathscr{L}$. The eigenvalue of the eigenvector $\mathbf{e}$ is
maximal, so the commutation relations lead to the following:
\begin{equation}\label{main_prop}
\mathscr{T} \left(c_{\overline{\omega}_{l_1},\lambda}^{\overline{\omega}_{l_1}}\right)
\mathbf{e}=0,\quad\lambda\ne
w_0\overline{\omega}_{l_1};\qquad\mathscr{T}_\mathcal{I}(z)\mathbf{e}=0, \quad
z\in\mathbb{C}[\mathfrak{p}^-]_{q,1}.
\end{equation}

Now we intend to introduce a Fock representation $\mathscr{T}_F$
of the $*$-algebra $\operatorname{Pol}(U)_q$ and a non-zero
intertwining operator between $\mathscr{T}_F$ and
$\mathscr{T}_\mathcal{I}$. Then the faithfulness of
$\mathscr{T}_\mathcal{I}$ will follow from the
faithfulness and irreducibility of $\mathscr{T}_F$.\\

The well-known equality $(S\otimes\mathrm{id})R=(\mathrm{id}\otimes S^{-1})R=R^{-1}$ (see
\cite[p. 326]{Drinf2}) leads to the invertibility of the linear map
$\check{R}_{\mathbb{C}[\mathfrak{p}^+]_q,\mathbb{C}[\mathfrak{p}^-]_q}$. This allows us
to turn from the expansion \eqref{expan1} to the expansion
$$
\operatorname{Pol}(U)_q= \bigoplus_{i,j=0}^\infty
\mathbb{C}[\mathfrak{p}^+]_{q,-j}\mathbb{C}[\mathfrak{p}^-]_{q,i}.
$$
Obviously, $\mathbb{C}[\mathfrak{p}^+]_{q,0}=
\mathbb{C}[\mathfrak{p}^-]_{q,0}=\mathbb{C}\cdot 1$. Hence for any
$f\in\operatorname{Pol}(U)_q$ there exists a unique decomposition similar to
\eqref{expan2}:
$$
f=\langle f\rangle\cdot 1+\sum_{(i,j)\ne(0,0)}f_{ij},\qquad
f_{ij}\in\mathbb{C}[\mathfrak{p}^+]_{q,-j}\mathbb{C}[\mathfrak{p}^-]_{q,i},
$$
where $\langle f\rangle$ is a linear functional on $\operatorname{Pol}(U)_q$.

Following section \ref{finite}, we introduce the
$\operatorname{Pol}(U)_q$-module $\mathscr{H}$ with a generator
$v_0$ and defining relations $fv_0=0$ for any
$f\in\mathbb{C}[\mathfrak{p}^-]_{q,1}$, i.e. $fv_0=\langle
f\rangle v_0$ for $f\in\mathbb{C}[\mathfrak{p}^-]_q$. Evidently,
$$
\mathbb{C}[\mathfrak{p}^+]_q\overset{\approx}{\to}\mathscr{H},\qquad
f\mapsto fv_0,
$$
which allows us to equip $\mathscr{H}$ with the sesquilinear form
\begin{equation}\label{h_f}
(f_1v_0,f_2v_0)=\langle f_2^\star f_1\rangle,\qquad
f_1,f_2\in\mathbb{C}[\mathfrak{p}^+]_q.
\end{equation}
It follows from \eqref{main_prop} that
$$
\langle
f\rangle=(\mathscr{T}_\mathcal{I}(f)\mathbf{e},\mathbf{e}),\qquad
f\in\operatorname{Pol}(U)_q.
$$
Hence \eqref{h_f} is a nonnegative Hermitian form. Its positivity can be proved using
methods from the function theory on compact quantum groups. Indeed, if $\langle f^\star
f\rangle=0$ and $f\in\mathbb{C}[\mathfrak{p}^+]_q$, then for $k\in\mathbb{N}$ large
enough, the element $\mathcal{I}(f)t^{\star k}$ belongs to $(\mathbb{C}[G]_q,\star)$ and
is annihilated by all its irreducible representations, see \cite[p. 121]{KorSoib}. Hence
the $L^2$-norm of $\mathcal{I}(f)t^{\star k}$ is equal to 0. The orthogonality relations
imply that $\mathcal{I}(f)t^{\star k}=0$, and Proposition \ref{hol_embed} leads to $f=0$.

Thus, $\mathscr{H}$ is a pre-Hilbert space and a $\operatorname{Pol}(U)_q$-module. Denote
by $\mathscr{T}_F$ the corresponding $*$-representation of $\operatorname{Pol}(U)_q$ in
$\mathscr{H}$; this is the mentioned Fock representation.

For obvious reasons, there exists a unique linear operator
$\mathscr{I}:\mathscr{H}\to\mathscr{L}$ which intertwines the
representations $\mathscr{T}_F$ and $\mathscr{T}_\mathcal{I}$, and such that
$\mathscr{I}:v_0\mapsto\mathbf{e}$.

The irreducibility of $\mathscr{T}_F$ can be proved similarly to the
irreducibility of $T_F$ in section \ref{finite}. Indeed, if
$\mathscr{H}'\ne 0$ is a common invariant subspace of
$\mathscr{T}_F(f)$, $f\in\operatorname{Pol}(U)_q$, then
$\mathscr{H}'$ has a non-zero vector $v$ such that
$$
\mathscr{T}_F(f)v=0,\qquad
f\in\bigoplus_{j=1}^\infty\mathbb{C}[\mathfrak{p}^-]_{q,j}.
$$
Hence $v$ is orthogonal to the subspace
$\left(\bigoplus\limits_{j=1}^\infty
\mathbb{C}[\mathfrak{p}^+]_{q,-j}\right)v_0$ and, therefore,
$\mathbb{C}v=\mathbb{C}v_0$, $\mathscr{H}'=\mathscr{H}$. The
irreducibility of $\mathscr{T}_F$ is thus proved.

We equip $\mathscr{H}$ with the grading
$$
\mathscr{H}=\bigoplus\limits_{j=0}^\infty\mathscr{H}_j,\qquad
\mathscr{H}_j=\mathbb{C}[\mathfrak{p}^+]_{q,-j}v_0.
$$
The subspaces $\mathscr{H}_j$ are pairwise orthogonal and finite
dimensional. As in section \ref{finite}, one can prove that for
any $j,k\in\mathbb{Z}_+$ the map
$$
\mathbb{C}[\mathfrak{p}^+]_{q,-j}\otimes\mathbb{C}[\mathfrak{p}^-]_{q,k}
\overset{m}{\to}
\mathbb{C}[\mathfrak{p}^+]_{q,-j}\mathbb{C}[\mathfrak{p}^-]_{q,k}\to
\operatorname{Hom}(\mathscr{H}_k,\mathscr{H}_j),\qquad
f\mapsto\mathscr{T}_F(f)|_{\mathscr{H}_k},
$$
is bijective (see Lemma \ref{ball_I_8.7}). Now use the positivity
of the Hermitian form $(\cdot,\cdot)$ in $\mathscr{H}$ to obtain the
faithfulness of $\mathscr{T}_F$ (see Lemma \ref{ball_I_8.8}).

The irreducibility and faithfulness of $\mathscr{T}_F$, hence the faithfulness of
$\mathscr{T}_\mathcal{I}$ and the injectivity of $\mathcal{I}$ are now proved.
\hfill $\square$

\bigskip

The next statement follows from the fact that $\mathbb{C}[G]_q$ is an integral domain.
\begin{corollary}\label{Pol-integrity} The algebra
$\operatorname{Pol}(\mathfrak{p}^-)_q$ is an integral domain.
\end{corollary}


\medskip

\subsection{\boldmath $U_q\mathfrak{k}$-invariant
polynomials}\label{inv_polynomials}

We introduce pairwise commuting $U_q\mathfrak{k}$-invariants $y_1,\ldots,y_r \in
\mathbb{C}[\mathfrak{p}^-]_q$ to be used extensively in the sequel. Here $r$ denotes the
rank of an irreducible bounded symmetric domain $\mathbb{D}$.

\begin{lemma}\label{spec_H_S} The set of all eigenvalues of the linear operator
$H_{\mathbb{S}}$ in $L(\overline{\omega}_{l_1})$ is the set
$$
 \{\,\overline{\omega}_{l_1}(H_\mathbb{S})\,,\,
\overline{\omega}_{l_1}(H_\mathbb{S})-2\,,\,\ldots\,,\,
-\overline{\omega}_{l_1}(H_\mathbb{S})+2\,,\,
-\overline{\omega}_{l_0}(H_\mathbb{S})\,\}.
$$
The number of eigenvalues is $r+1$, with $r$ being the rank of the
bounded symmetric domain $\mathbb{D}$.
\end{lemma}

{\bf Proof.} The highest weight of the $U_q\mathfrak{g}$-module
$L(\overline{\omega}_{l_1})$ equals $\overline{\omega}_{l_1}$, the lowest one equals to
$-\overline{\omega}_{l_0}$. Hence the first statement follows from the simplicity of the
$U_q\mathfrak{g}$-module $L(\overline{\omega}_{l_1})$ and the definition of
$H_\mathbb{S}\in\mathfrak{h}$. The second statement is equivalent to
\begin{equation}\label{2_r}
\overline{\omega}_{l_1}(H_\mathbb{S})+\overline{\omega}_{l_0}(H_\mathbb{S})
=2r.
\end{equation}
The last equation is obtainable from the Cartan list of irreducible
bounded symmetric domains \cite{Helg}. Indeed, the ranks of these
domains are known from \cite[Sec. 9.4.4]{Helg}, and the summands
in the left hand side of \eqref{2_r} are just the coefficients of
$\alpha_{l_0}$ in the expansions of the fundamental weights
$\overline{\omega}_{l_1}$, $\overline{\omega}_{l_0}$ into sums of
simple roots of $\mathfrak{g}$. These coefficients are listed in
\cite{Bou4-6}. \hfill $\square$

\medskip

Recall the notation introduced in section \ref{true_bases}
$$
v_{\lambda,j}^{\overline{\omega}_{l_1}},\qquad
c_{\overline{\omega}_{l_1};\lambda,j}^{\overline{\omega}_{l_1}}
$$
for the elements of a special basis of
$L(\overline{\omega}_{l_1})$ and the matrix elements of the
corresponding fundamental representation in this basis.

\begin{lemma}\label{z_lambda_j}
For any matrix element
$c_{\overline{\omega}_{l_1};\lambda,j}^{\overline{\omega}_{l_1}}$
there is a unique element $z_{\lambda,j}$ of
$\mathbb{C}[\mathfrak{p}^-]_q$ such that
$$
\mathcal{I}(z_{\lambda,j})=t^{-1}\cdot
c_{\overline{\omega}_{l_1};\lambda,j}^{\overline{\omega}_{l_1}}.
$$
\end{lemma}

{\bf Proof.} The uniqueness follows from the injectivity of
$\mathcal{I}$. To prove the existence of $z_{\lambda,j}$, consider
the subspace $t \cdot
\mathcal{I}\mathbb{C}[\mathfrak{p}^-]_q\subset\mathbb{C}[G]_{q,x}$.
It contains $t$ and is a $U_q\mathfrak{b}^+$-submodule of the
$U_q\mathfrak{b}^+$-module $\mathbb{C}[G]_{q,x}$ since for any
$j=1,2,\ldots,l$
\begin{gather*}
K_j^{\pm 1}(t\mathcal{I}(f))=(K_j^{\pm 1}t)\mathcal{I}(K_j^{\pm 1}f),
\\ E_j(t\mathcal{I}(f))=(E_jt)\mathcal{I}f+(K_jt)\mathcal{I}(E_jf),\qquad
f\in\mathbb{C}[\mathfrak{p}^-]_q,
\end{gather*}
and $K_j^{\pm 1}t$, $E_jt \in t\cdot \mathcal{I}\mathbb{C}[\mathfrak{p}^-]_q$.
Hence,
$$
t \cdot \mathcal{I}\mathbb{C}[\mathfrak{p}^-]_q\supset
U_q\mathfrak{b}^+\cdot t,
$$
and it now just remains to observe that
$c_{\overline{\omega}_{l_1};\lambda,j}^{\overline{\omega}_{l_1}}
\in U_q\mathfrak{b}^+\cdot t$. \hfill $\square$

\medskip

For short, we use the notation $z_\lambda$ instead of $z_{\lambda,1}$ in case the weight
subspace $L(\overline{\omega}_{l_1})_\lambda$ is one dimensional.

\medskip

Lemmas \ref{spec_H_S}, \ref{z_lambda_j} make it possible to introduce the
elements
\begin{equation}\label{def_y_i}
y_i=\sum_{\{(\lambda,j)|\:
\lambda(H_\mathbb{S})=-\overline{\omega}_{l_0}(H_\mathbb{S})+2i\}}
z_{\lambda,j}\cdot z_{\lambda,j}^*,\qquad i=1,\ldots,r.
\end{equation}
Obviously, $y_i=y_i^*$ for any $i=1,2,\ldots,r$. Put $y_0=1$.

\begin{proposition}\label{commut_y_i}
1. The elements $y_1,y_2,\ldots,y_r$ of
$\operatorname{Pol}(\mathfrak{p}^-)_q$ are
$U_q\mathfrak{k}$-invariant and do not depend on the choice of the
orthonormal basis
$\left\{v_{\lambda,j}^{\overline{\omega}_{l_1}}\right\}$ in
$L(\overline{\omega}_{l_1})$.
\begin{flalign*}
\it{2}. && y_iy_j=y_jy_i,\qquad i,j=1,2,\ldots,r.&&
\end{flalign*}
\end{proposition}

{\bf Proof.} The first statement follows from the facts that
$\operatorname{Pol}(\mathfrak{p}^-)_q$ is a
$(U_q\mathfrak{k},*)$-module algebra, the Hermitian form in
$L(\overline{\omega}_{l_1})$ is $U_q\mathfrak{k}$-invariant, and
the basis $\left\{v_{\lambda,j}^{\overline{\omega}_{l_1}}\right\}$
is orthonormal.

The second statement follows from the $U_q\mathfrak{k}$-invariance of
$y_i$, together with the following result obtained by D.~Shklyarov.

\begin{lemma}\label{k-invariants}
The subalgebra of $U_q\mathfrak{k}$-invariants in
$\operatorname{Pol}(\mathfrak{p}^-)_q$ is commutative.
\end{lemma}

{\bf Proof.} Consider the $U_q\mathfrak{k}$-submodule
$\mathbb{C}[\mathfrak{p}^-]_qf_0\subset D(\mathbb{D})_q$. Using the
Hua-Schmid theorem (see \cite[p. 73]{John}, \cite[p. 443]{Takeuchi}), one
can prove that, as a $U_q\mathfrak{k}$-module,
$$
\mathbb{C}[\mathfrak{p}^-]_qf_0\approx
\bigoplus_{\lambda\in\mathcal{A}_+}L(\mathfrak{k},\lambda),
$$
where $\mathcal{A}_+=\bigoplus_{j=1}^r \mathbb Z_+ (\sum_{i=1}^j \gamma_j)\subset
P_+^\mathbb{S}$ with $\gamma_1,...,\gamma_r$ being positive noncompact strongly
orthogonal roots. Hence any endomorphism of the $U_q\mathfrak{k}$-module
$\mathbb{C}[\mathfrak{p}^-]_q$ acts as a scalar multiplication in any
$U_q\mathfrak{k}$-isotypic components of $\mathbb{C}[\mathfrak{p}^-]_q$. Therefore the
algebra of endomorphisms of the $U_q\mathfrak{k}$-module $\mathbb{C}[\mathfrak{p}^-]_q$
is commutative.

In particular, the operators of multiplication by $\varphi$, $\psi$
$$f\mapsto \varphi f,\qquad f\mapsto\psi f,$$
commute if $\varphi,\psi\in\operatorname{Pol}(\mathfrak{p}^-)_q$ are
$U_q\mathfrak{k}$-invariant. This means that the operators of the Fock representation
$T_F(\varphi)$, $T_F(\psi)$ of $\operatorname{Pol}(\mathfrak{p}^-)_q$ in the pre-Hilbert
space $\mathcal{H}\cong\mathbb{C}[\mathfrak{p}^-]_qf_0$ commute. The result then follows
from Lemma \ref{ball_I_8.8} on the faithfulness of the Fock representation. \hfill
$\square$

\begin{remark}\label{Fourier}
In the category of $U_q\mathfrak{k}$-modules,
$$
\mathcal{H}=\bigoplus_{\lambda\in\mathcal{A}_+}\mathcal{H}_\lambda,\qquad
\mathcal{H}_\lambda \approx L(\mathfrak{k},\lambda),
$$
and any $U_q\mathfrak{k}$-invariant
$\psi\in\operatorname{Pol}(\mathfrak{p}^-)_q$ is determined by its `Fourier
coefficients' $\widehat{\psi}(\lambda)=T_F(\psi)|_{\mathcal{H}_\lambda}$
(see the proof of Lemma \ref{k-invariants}).
\end{remark}


\medskip

\subsection{Representations of the $*$-algebra \boldmath
$\operatorname{Pol}(\mathfrak{p}^-)_q$.} \label{reps_pol}

In this section we finish the proof of the existence and uniqueness of
a faithful irreducible $*$-representation of
$\operatorname{Pol}(\mathfrak{p}^-)_q$ by bounded operators in a
Hilbert space.

Use $y_1,y_2,\ldots, y_r$ introduced in the last section.

\begin{proposition}\label{def_y}
There exists a unique element $y$ in $\operatorname{Pol}(\mathfrak{p}^-)_q$
such that $\mathcal{I}y=x^{-1}$. It is $U_q\mathfrak{k}$-invariant and is
given by
\begin{equation}\label{explicit_y}
y=1+\sum_{i=1}^r(-1)^iy_i.
\end{equation}
\end{proposition}

{\bf Proof.} The uniqueness of $y$ is obvious. We now show that its existence and
\eqref{explicit_y} follow from the orthogonality relations. Indeed, the orthogonality
relations imply that
\begin{equation}\label{unitary_new}
\sum_{(\lambda,j)}
c_{\overline{\omega}_{l_1};\lambda,j}^{\overline{\omega}_{l_1}}\cdot
\left(c_{\overline{\omega}_{l_1};\lambda,j}^{\overline{\omega}_{l_1}}\right)
^\star=1.
\end{equation}
But
\begin{equation}\label{star_aster_new}
\left(c_{\overline{\omega}_{l_1};\lambda,j}^{\overline{\omega}_{l_1}}\right)
^\star=(-1)^{\frac{\lambda(H_\mathcal{S})+\lambda(\overline{\omega}_{l_0})}2}
\left(c_{\overline{\omega}_{l_1};\lambda,j}^{\overline{\omega}_{l_1}}\right)
^*
\end{equation}
since
$$
t^\star=t^*,\qquad F_j^\star=
\begin{cases}
F_j^*, & j\ne l_0,
\\-F_j^*, & j=l_0,
\end{cases}
$$
and
$$
(\xi f)^\star=(S(\xi))^\star f^\star,\qquad(\xi f)^*=(S(\xi))^*f^*
$$
for any $\xi\in U_q\mathfrak{g}$, $f\in\mathbb{C}[G]_q$. It follows from
\eqref{unitary_new}, \eqref{star_aster_new} that
$$
t\mathcal{I}\left(\sum_{i=0}^r(-1)^iy_i\right)t^*=1.
$$
Now it just remains to multiply both sides of the last equality by
$t^{-1}$ and $(t^*)^{-1}$.

\eqref{explicit_y} and Proposition \ref{commut_y_i} imply the
$U_q\mathfrak{k}$-invariance of $y$. \hfill $\square$

\medskip

\begin{corollary}\label{quasicom_y}
For any $z\in\mathbb{C}[\mathfrak{p}^-]_{q,1}$ one has
\begin{equation}\label{q_y}
z\, y=q_{l_0}^{-2}\,y\, z,\qquad z^*\, y=q_{l_0}^2\, y\, z^*.
\end{equation}
\end{corollary}

{\bf Proof.} In the special case $z=z_{\rm low}$ the commutation relations \eqref{q_y}
are easily deducible from the equality $\mathcal{I}y=x^{-1}$ and the commutation
relations for the matrix elements in a special basis, together with Corollary
\ref{quasicom_x}. The general case reduces to the special one via the action of
$U_q\mathfrak{k}$. \hfill $\square$

\bigskip

Recall the notation $T_F$ for the Fock representation of
$\operatorname{Pol}(\mathfrak{p}^-)_q$ in the pre-Hilbert space
$\mathcal{H}$ and its completion $\overline{\mathcal{H}}$. As is proved in
\cite{SSV4}, the operators $T_F(f)$ are bounded, hence $T_F$ admits an
extension by continuity to a representation $\overline{T}_F$ in the Hilbert
space $\overline{\mathcal{H}}$.

\begin{proposition}\label{irreps_Pol}
1. The representation $\overline{T}_F$ is a faithful irreducible $*$-representation
of $\operatorname{Pol}(\mathfrak{p}^-)_q$ by bounded operators in a Hilbert
space.

2. A representation with these properties is unique up to a unitary equivalence.
\end{proposition}

{\bf Proof.} 1. The faithfulness of $\overline{T}_F$ follows from the
faithfulness of $T_F$ that is proved in section \ref{finite}. Now let us
turn to the proof of irreducibility.

It follows from \eqref{def_y_i} that the self-adjoint linear
operator $\overline{T}_F(y)$ is compact and its spectrum is a
closure of a geometric series. More precisely,
\begin{equation}\label{spec_y}
\overline{T}_F(y)|_{\mathcal{H}_k}=q_{l_0}^{2k}\cdot 1,\qquad k\in\mathbb{Z}_+,
\end{equation}
with $\mathcal{H}_k=\mathbb{C}[\mathfrak{p}^-]_{q,k}\cdot v_0$.
Let $\mathfrak{A}$ be the closure of the algebra of operators
$\left\{\left.\overline{T}_F(f)\right|\:
f\in\operatorname{Pol}(\mathfrak{p}^-)_q\right\}$ with respect to
the operator norm. It follows from \eqref{spec_y} that the algebra
$\mathfrak{A}$ contains the orthogonal projections onto all the
subspaces $\mathcal{H}_k$. These subspaces are pairwise
orthogonal, finite dimensional, and their sum is dense in
$\overline{\mathcal{H}}$. Hence Lemma \ref{ball_I_8.7} implies
that $\mathfrak{A}$ contains all compact linear operators in
$\overline{\mathcal{H}}$. This leads to the irreducibility of
$\overline{T}_F$.

Obviously, $\overline{T}_F$ is a $*$-representation.

\medskip

2. Let $T'$ be a faithful irreducible $*$-representation of
$\operatorname{Pol}(\mathfrak{p}^-)_q$ by bounded linear operators in a
Hilbert space. We now show that $T'$ and $\overline{T}_F$ are unitary equivalent.

The faithfulness of $T'$ implies that $T'(y)\neq 0$. This allows us to use the
commutation relations \eqref{q_y} to prove the fact that the spectrum of the self-adjoint
linear operator $T'(y)$ is the closure of a geometric series.

Let $v'$ be a normalized eigenvector with eigenvalue of maximum
modulus. Then \eqref{q_y} implies $T'(z)^*v'=0$ for any
$z\in\mathbb{C}[\mathfrak{p}^-]_{q,1}$. Hence
$(T'(f)v',v')=(\overline{T}_F(f)v_0,v_0)$ for any
$f\in\operatorname{Pol}(\mathfrak{p}^-)_q$. Therefore the map
$v_0\mapsto v'$ admits an extension to an isometric linear
operator that intertwines $\overline{T}_F$ and $T'$. The operator
is surjective since $T'$ is irreducible.\hfill $\square$

\begin{remark}
In the proof of uniqueness we only use the fact that $T'(y)\ne 0$,
not the faithfulness of $T'$.
This observation has an interesting application. Let $J=J^*$ be a
non-zero two-sided ideal of $\mathrm{Pol}(\mathfrak{p}^-)_q$. If
the factor algebra $\mathrm{Pol}(\mathfrak{p}^-)_q/J$ has
irreducible $*$-representations that separate its points, then
$y\in J$. In the special case of the quantum matrix ball and the
defining ideal $J$ of its Shilov boundary \cite[p.
381-383]{Vak01}, this implies the inclusion $(y) \subset J$, which
means that the Shilov boundary belongs to the topological boundary
of the ball. Here $(y)$ stands for the two-sided ideal generated
by $y$.
\end{remark}

\bigskip

The multiplicative system $y^{\mathbb{Z}_+}$ is an Ore set of the
integral domain $\operatorname{Pol}(\mathfrak{p}^-)_q$, see
\eqref{q_y} and Corollary \ref{Pol-integrity}. Consider the
localization $\operatorname{Pol}(\mathfrak{p}^-)_{q,y}$ of the
algebra $\operatorname{Pol}(\mathfrak{p}^-)_q$ with respect to the
set $y^{\mathbb{Z}_+}$. Evidently,
$\operatorname{Pol}(\mathfrak{p}^-)_q\hookrightarrow
\operatorname{Pol}(\mathfrak{p}^-)_{q,y}$. Moreover, the natural
extension of $T_F$ from $\operatorname{Pol}(\mathfrak{p}^-)_q$ to
$\operatorname{Pol}(\mathfrak{p}^-)_{q,y}$ is a faithful
representation of $\operatorname{Pol}(\mathfrak{p}^-)_{q,y}$ in
$\mathcal{H}$. Similarly to the proof of Theorem \ref{c_e}, one
has

\begin{proposition}\label{inj}
The natural extension of the canonical embedding
$\operatorname{Pol}(\mathfrak{p}^-)_q\hookrightarrow
\mathbb{C}[\mathbb{X}_\mathbb{S}]_{q,x}$ to the localization
$\operatorname{Pol}(\mathfrak{p}^-)_{q,y}$ is an embedding of algebras
\begin{equation}\label{ext_inj}
\operatorname{Pol}(\mathfrak{p}^-)_{q,y}\hookrightarrow
 \mathbb{C}[\mathbb{X}_\mathbb{S}]_{q,x}.
\end{equation}
\end{proposition}

\begin{remark}
As in the previous sections, one can prove that the
$(U_q\mathfrak{g},*)$-module algebra structure is canonically
extendable from $\operatorname{Pol}(\mathfrak{p}^-)_q$ to
$\operatorname{Pol}(\mathfrak{p}^-)_{q,y}$. Also, \eqref{ext_inj}
is a morphism of $(U_q\mathfrak{g},*)$-module algebras (we use the
fact that the elements $\xi(fy^n)y^{-n}$, $\xi\in
U_q\mathfrak{g}$, $f\in\operatorname{Pol}(\mathfrak{p}^-)_q$, are
values of the corresponding Laurent polynomials on a geometric
series).
\end{remark}

\bigskip

We maintain the notation $\mathcal{I}$ for {\it the canonical embedding}
\eqref{ext_inj} and now describe the image of
$\operatorname{Pol}(\mathfrak{p}^-)_{q,y}$ under it.

Let $K$ be the connected affine algebraic subgroup of $G$ with Lie algebra
$\mathfrak{k}$, see section \ref{finite}, $w_0$ the longest element of $W$, and
\begin{equation*}
K_1=w_0 K w_0^{-1}.
\end{equation*}
Obviously, if $a,b \in G$ then $w_0\cdot Ka=w_0\cdot Kb$ is
equivalent to $K_1\,(w_0a)=K_1\,(w_0b)$.

Solutions of the system of equations
$$ L_{\rm reg}(E_i)f=L_{\rm reg}(F_i)f=0,\qquad i\neq l_1,$$
$$ L_{\rm reg}(K_j^{\pm 1})f=f,\qquad j=1,2,\ldots,l$$
form a subalgebra in the $U_q\mathfrak{g}$-module algebra $\mathbb{C}[G]_q$. It can be
treated as a quantum analog of the algebra $\mathbb{C}[K_1\backslash G]$ of regular
functions on the affine algebraic variety $K_1\backslash G$; see Proposition
\ref{from_Jo-R}. Denote by $\mathbb{C}[K_1\backslash G]_q$ the latter
$U_q\mathfrak{g}$-module algebra.

The algebra $\mathbb{C}[K_1\backslash G]_q$ is generated by
$c_{\overline{\omega}_{l_1};\lambda,j}^{\overline{\omega}_{l_1}}
(c_{\overline{\omega}_{l_1};\mu,k}^{\overline{\omega}_{l_1}})^\star$. Hence it can be
defined as the smallest $U_q\mathfrak{g}$-module subalgebra that contains $x$. This
definition is used implicitly in the sequel.

\medskip

In section \ref{canon_embed} the algebra $\mathbb{C}[\mathbb{X}_\mathbb{S}]_{q}$ is
equipped with a grading. We now extend it to the localization
$\mathbb{C}[\mathbb{X}_\mathbb{S}]_{q,x}$.

\begin{lemma}\label{gradation}
There exists a unique $\mathbb{Z}$-grading of $\mathbb{C}[\mathbb{X}_\mathbb{S}]_{q,x}$
such that $\deg\left(c_{\overline{\omega}_{l_1};\lambda,j}
^{\overline{\omega}_{l_1}}\right)=1$,\; $\deg\left(c_{-\overline{\omega}_{l_1};\mu,k}
^{\overline{\omega}_{l_0}}\right)=-1.$
\end{lemma}

{\bf Proof.} The uniqueness is obvious. A construction of a grading with the required
properties can be done as in section \ref{w0_G_x}. Indeed,
$$
L_\mathrm{reg}(H_{l_1})c_{\overline{\omega}_{l_1};\lambda,j}
^{\overline{\omega}_{l_1}}=c_{\overline{\omega}_{l_1};\lambda,j}
^{\overline{\omega}_{l_1}},\qquad
L_\mathrm{reg}(H_{l_1})c_{-\overline{\omega}_{l_1};\mu,k}
^{\overline{\omega}_{l_0}}=-c_{-\overline{\omega}_{l_1};\mu,k}
^{\overline{\omega}_{l_0}},
$$
and the differentiation $L_\mathrm{reg}(H_{l_1})$ on
$\mathbb{C}[\mathbb{X}_\mathbb{S}]_{q}$ admits an extension to a
differentiation on $\mathbb{C}[\mathbb{X}_\mathbb{S}]_{q,x}$. It
remains to set $\deg f=j$ for any $f
\in\mathbb{C}[\mathbb{X}_\mathbb{S}]_{q,x}$ such that
$L_\mathrm{reg}(H_{l_1})f=jf$.\hfill $\square$

\medskip

Note that the localization $\mathbb{C}[K_1\backslash G]_{q,x}$
of $\mathbb{C}[K_1\backslash G]_q$ with respect to the
multiplicative system $x^{\mathbb{Z}_+}$ is naturally embedded
into $\mathbb{C}[\mathbb{X}_\mathbb{S}]_{q,x}$.

\begin{proposition}\label{image_emb}\ \ \
$\mathcal{I}\operatorname{Pol}(\mathfrak{p}^-)_{q,y}=
\mathbb{C}[K_1\backslash G]_{q,x}$.
\end{proposition}

{\bf Proof.} The embedding of $U_q\mathfrak{g}$-module algebras
$\mathcal{I}\operatorname{Pol}(\mathfrak{p}^-)_{q,y}\subset
\mathbb{C}[K_1\backslash G]_{q,x}$ is obvious. To prove the
opposite inclusion, consider the subspace
$$
\left(\bigoplus_{j=1}^\infty\,t^j\cdot\mathcal{I}
\operatorname{Pol}(\mathfrak{p}^-)_{q,y}\right)\;\bigoplus\;
\mathcal{I}\operatorname{Pol}(\mathfrak{p}^-)_{q,y}\;\bigoplus\;
\left(\bigoplus_{j=1}^\infty\,
\mathcal{I}\operatorname{Pol}(\mathfrak{p}^-)_{q,y}\cdot
t^{*j}\right).
$$
This is a $U_q\mathfrak{g}$-module subalgebra; it contains $t$, $t^*$,
$x^{-1}$. Hence it coincides with
$\mathbb{C}[\mathbb{X}_\mathbb{S}]_{q,x}$. Therefore
$\mathcal{I}\operatorname{Pol}(\mathfrak{p}^-)_{q,y}$ contains all elements
of degree zero from $\mathbb{C}[\mathbb{X}_\mathbb{S}]_{q,x}$.
It remains now just to use the fact that non-zero elements of
$\mathbb{C}[K_1\backslash G]_{q,x}$ have zero degree. \hfill $\square$


\medskip

\subsection{Appendix. A cone over $\partial \mathbb{D}$.}

Consider the two-sided ideal $(y)$ of the algebra
$\operatorname{Pol}(\mathfrak{p}^-)_q$ generated by $y=y^*$, and
the corresponding factor-algebra
$$ \mathbb{C}[\partial \mathbb{D}]_q\;=\;
\operatorname{Pol}(\mathfrak{p}^-)_q/(y).$$

It is a $q$-analog for the algebra of regular functions on the
affine algebraic variety $\partial \mathbb{D}$ which is the
topological boundary of $\mathbb{D}$. It is important to note that
the two-sided ideal $(y)$ is a $U_q\mathfrak{g}$-module subalgebra
of the $U_q\mathfrak{g}$-module algebra
$\operatorname{Pol}(\mathfrak{p}^-)_q$. Indeed,
$$ E_jy=F_jy=0,\quad j\ne l_0,\qquad\qquad K_i^{\pm 1}y=y,\quad i=1,2,\ldots,l,$$
and it remains to prove that $E_{l_0}y\in (y)$, $F_{l_0}y\in (y)$.
The last statement can be verified easily by using the equality
$\mathcal{I}\,y=t^{-1}(t^*)^{-1}$. Hence $\mathbb{C}[\partial
\mathbb{D}]_q$ is a $(U_q\mathfrak{g},*)$-module algebra.

Consider the auxiliary $*$-subalgebra of
$\mathbb{C}[\mathbb{X}]_{q,x}$ generated by the set $\mathcal{I}
\mathbb{C}[\mathfrak{p}^-]_{q,1}$ and $t$. It is a
$(U_q\mathfrak{g},*)$-module subalgebra since
$$ E_jt=F_jt=(K_j^{\pm 1}-1)t=0,\qquad j\ne l_0,
$$
$$ F_{l_0}t=(K_{l_0}^{\pm 1}-1)t=0,\qquad E_{l_0}t=q_{l_0}^{-1/2}\,t\,
(\mathcal{I}z_{\mathrm{low}}).
$$

One can describe the latter $*$-algebra by its generators $t$,
$z\in \mathbb{C}[\mathfrak{p}^-]_{q,1}$ and the relations
\begin{equation}\label{4.07-1}
tt^*=t^*t,\qquad tt^*y=ytt^*=1,
\end{equation}
\begin{equation}\label{4.07-2}
tz=q_{l_0}^{-1}zt,\quad t^*z=q_{l_0}^{-1}zt,\qquad\qquad z\in
\mathbb{C}[\mathfrak{p}^-]_{q,1}.
 \end{equation}
Moreover, the set of relations should include the defining
relations of the $*$-algebra
$\operatorname{Pol}(\mathfrak{p}^-)_q$.

To obtain \eqref{4.07-1}, \eqref{4.07-2}, one can identify an element $z\in
\mathbb{C}[\mathfrak{p}^-]_{q,1}$ with its image under the canonical embedding
$\mathcal{I}$.

The action of $U_q\mathfrak{g}$ is given as before
$$ E_jt=F_jt=(K_j^{\pm 1}-1)t=0,\qquad j\ne l_0,$$
$$
F_{l_0}t=(K_{l_0}^{\pm 1}-1)t=0,\qquad E_{l_0}t=q_{l_0}^{-1/2}\,t\,
(\mathcal{I}z_{\mathrm{low}}).
$$


\end{document}